\documentclass[11pt, reqno]{amsart}

\usepackage{amssymb,amsmath,amscd,graphicx,
latexsym,amsthm}
\usepackage{amssymb,latexsym,amsmath,amscd,graphicx}
\usepackage[mathscr]{eucal}

\usepackage{scalerel,stackengine}
\DeclareMathSymbol{\shortminus}{\mathbin}{AMSa}{"39}
\usepackage{multicol}
\usepackage{cite}

\usepackage{xcolor}

\stackMath
\newcommand\reallywidehat[1]{%
\savestack{\tmpbox}{\stretchto{%
  \scaleto{%
    \scalerel*[\widthof{\ensuremath{#1}}]{\kern.1pt\mathchar"0362\kern.1pt}%
    {\rule{0ex}{\textheight}}%WIDTH-LIMITED CIRCUMFLEX
  }{\textheight}% 
}{2.4ex}}%
\stackon[-6.9pt]{#1}{\tmpbox}%
}
 \input xy
 \xyoption{all}
\def\OO{{\mathcal O}}

\def\F{\mathcal{F}}
\def\E{\mathcal{E}}

\def\I{\mathcal{I}}

\def\cP{\mathcal{P}}
\def\Pic0{{\rm Pic}^0}

\def\Aut0{{\rm Aut}^0}

\def\T{\mathcal T} 
\def\l{{\underline l}}

\def\*{{\underline *}}
\def\utheta{{\underline \theta}}

\def\rk{\mathrm{rk}}

\def\udelta{{\underline \delta}}

\theoremstyle{plain}

\newtheorem{theorem}{Theorem}[subsection]
\newtheorem{theoremalpha}{Theorem}
\newtheorem{corollaryalpha}[theoremalpha]{Corollary}

\newtheorem{proposition/example}[theorem]{Proposition/Example}
\newtheorem{proposition}[theorem]{Proposition}

\newtheorem{lemma}[theorem]{Lemma}

\theoremstyle{definition}

\newtheorem{remark}[theorem]{Remark}

\newtheorem{conjecture/question}[theorem]{Conjecture/Question}

\newtheorem{remark/definition}[theorem]{Remark/Definition}
\newtheorem{notation/assumptions}[theorem]{Assumptions/Notation}
\newtheorem{setting/notation}[theorem]{Setting and Notation}

\numberwithin{equation}{section}

 \theoremstyle{remark}

\begin{document}
\title[Abelian varieties analogs of two results about algebraic curves]{Abelian varieties analogs of two results about algebraic curves}
 \author[N.Alvarado, G.Pareschi]{Nelson Alvarado and Giuseppe Pareschi}

\address{ Departamento de Matemáticas, Facultad de Ciencias, Universidad de Chile, Las Palmeras 3425,  Santiago\\Chile}
\email{nelson.alvarado@ug.uchile.cl}
\address{Department of Mathematics,
              University of Rome Tor Vergata\\Italy}
\email{pareschi@mat.uniroma2.it}
 \thanks{Both authors were partially supported by  the MIUR Excellence Department Project MatMod@TOV awarded to the Department of Mathematics of the University of Rome Tor Vergata. GP was also partially supported by the PRIN 2022 "Moduli spaces and Birational Geometry" and is a member of GNSAGA - INDAM}

%\subjclass[2010]{14K25; 32G20}
\begin{abstract} We characterize decomposable principally polarized abelian varieties of the form $E\times B$, with $E$ an elliptic curve, in two different ways,  which are, surprisingly, completely analogous to  classical results of curve theory concerning hyperelliptic curves. The first one  is by the failure of a normal generation property, namely the generation in degree zero of a certain graded module over the symmetric algebra over $H^0(2\Theta)$. 
 This appears to be the first result of this type in the realm of p.p.a.v.'s. The second characterization is by the failure of surjectivity of second order gaussian maps associated to line bundles corresponding to $6\Theta$, or, equivalently, by the fact that at some point, the line bundle corresponding to $3\Theta$ fails to separate $2$-jets.  
We also show that this last result is equivalent to an effective version of a  theorem of Nakamaye characterizing the above decomposable abelian varieties as those computing the minimal Seshadri constant. Finally we propose some conjectural generalizations relating $p$-jets separation thresholds, higher gaussian maps sujectivity thresholds, and Seshadri constants.
\end{abstract}

\maketitle 

\section{Introduction}\label{sec:intro}

In curve theory it is possible to characterize curves of genus $\ge 2$ with special Brill-Noether properties 
 by means of algebraic properties of graded section modules over the symmetric algebra associated to line bundles of suitably high degree. We allude to properties such as normal generation, normal presentation, linearity up to a certain stage of the minimal free resolution, i.e. M. Green's $N_p$ condition.  Another (less known) way of characterizing curves with special geometry is via the non-surjectivity of certain higher gaussian maps. An interesting feature of this sort of results is that they can be stated in terms of the failure of surjectivity of certain uniformly defined linear maps, thus providing equations in the moduli space of curves  for the loci in question.

%For several reasons, in higher dimension there is almost nothing like that, although an observation of Mukai sees some of the  fundamental theorems about normal generation and syzygies of curves  as Fujita-type results, and as such suited  to higher dimensional generalizations. Concerning principally polarized abelian varieties (the focus of this paper), the Fuijita program (originated from theorems of Koizumi and Mumford, and conjectured by Lazarsfeld) has been succesfully accomplished more than two decades ago: there are optimal theorems on normal generation, normal presentation, vanishing of Koszul cohomology groups, surjectivity of higher Wahl gaussian maps which are perfectly analogous to some theorems of curve theory. Despite this, with very few exceptions, there are no  characterizations of "special"  p.p.a.v.'s. via this sort of properties.  To begin with, it is not clear  what  geometric properties one should be able to distinguish in this way XXX.\footnote{
In this note we show two results suggesting that
 there is a chance that similar things can be done in the realm of 
 principally polarized abelian varieties.
 To put the main results into perspective, recall that, among smooth curves of fixed genus $g\ge 2$,   hyperelliptic curves are the most special ones in many respects (e.g. Brill-Noether theory). As such, they admit  many different characterizations. 
We will show that decomposable $g$-dimensional p.p.a.v.'s $E\times B$, with $\dim E=1$ seem to play the same role among p.p.a.v's, as they can be characterized in the same way. 
%%%%%%%%%%%%%%%%
\subsection{First characterization. }   Given a polarized projective variety $(X,L)$, with $L$ base point free, and a coherent sheaf $F$ on $X$, we denote
\begin{equation}\label{module}\mathcal R_{X,L}(F):=\bigoplus_{n\ge 0}H^0(X,F\otimes L^{\otimes n})
\end{equation}
which is a graded module over $\mathcal S_X(L)$, the symmetric algebra over $H^0(X,L)$. 
Our first result is

\begin{theoremalpha}\label{th:first} Let $(A,\OO_A(\Theta))$ be a principally polarized abelian variety. 
The following are equivalent:
\begin{enumerate}
\item $A$ is isomorphic, as polarized variety, to $E\times B$, where $\dim E=1$.
\item For every  $x\in A$ there exists  $\alpha\in \Pic0 A$ such that the graded module
\begin{equation}\label{eq:m(x,alpha)}
\mathcal R_{A,\OO_A(2\Theta)}(\I_x(3\Theta)\otimes P_\alpha),
\end{equation}
is not generated in degree $0$ over the symmetric algebra $\mathcal S_X(\OO_A(2\Theta))$, where $\I_{x}$ is the ideal sheaf of the point $x\in A$, and  $P_{\alpha}$ is the line bundle on $A$ parametrized by $\alpha.$ 
\item There exists  $x\in A$ and $\alpha\in \Pic0 A$ such that the graded module \emph{(\ref{eq:m(x,alpha)})} is not generated in degree $0$.
\end{enumerate}
\end{theoremalpha}

A more manageable version of Theorem \ref{th:first}, showing that the conditions (2) and (3) are equivalent to the non-surjectivity of a single uniformly defined linear map, is given by the following
\begin{corollaryalpha}\label{cor:B} Let $\Theta$ by a symmetric divisor representing the principal polarization, and let $e\in A$ be the neutral element.  Then $A\cong E\times B$, as in Theorem \ref{th:first}, if and only if the multiplication map of global sections
\begin{equation}\label{eq:0}
H^0(A,\OO_A(2\Theta))\otimes H^0(A,\I_e(3\Theta))\rightarrow H^0(A,\I_e(5\Theta))
\end{equation}
is not surjective.
\end{corollaryalpha}

The motivation of Theorem \ref{th:first} stems from classical results on curves and abelian varieties. 

\noindent\emph{Theorems of Castelnuovo, Mumford, Green-Lazarsfeld on curves and of Koizumi on  abelian varieties. } Let $C$ be a smooth projective curve of genus $g$. A result of Mumford (sligthly improving an earlier projective normality result of Castelnuovo) asserts that, as soon as $L$ and $M$ are line bundles on $C$ of degree $\ge 2g$ such that $\deg L+\deg M\ge 4g+1$ then the module 
$
\mathcal R_{C,L}(M)
$
 is generated in degree $0$ (\cite[Theorem 6]{mumford-quad}). According to Mukai, it is conceptually useful to write the line bundles $L$  and $M$ \`a la Fujita, i.e. $L=K_C(E)$, $M=K_C(F)$ where $K_C$ is the canonical bundle, so that the hypothesis of Mumford's  theorem is $\deg E,\deg F\ge 2$ and $\deg E+\deg F\ge 5$. 

In many respects abelian varieties  behave similarly to curves, and in fact there is a well known analog of Mumford's theorem for abelian varieties, due to Koizumi. As above, given an abelian variety $A$,  for $\alpha \in\Pic0 A$ we denote $P_\alpha$ the corresponding line bundle on $A$.  Koizumi's theorem (\cite[Theorem 4.6]{koizumi}), building on previous work of Mumford's (\cite{mumford-edav}),  states that if $(A,L)$ is a polarized abelian variety then the graded modules
\begin{equation}\label{eq:koizumi}
\mathcal R_{A,L^{\otimes n}}(L^{\otimes m}\otimes P_\alpha)
\end{equation} 
are generated in degree $0$ for all $\alpha\in \Pic0 A$ as soon as $n,m\ge 2$ and $n+m\ge 5$. 

 The theorems of Castelnuovo and Mumford (and, subsequently, Green's generalization to syzygies, \cite[Theorem 4.a.1]{green1}) have been the starting points of an important and far-reaching research direction in curve theory, consisting in characterizing  curves of given genus $g\ge 2$ with special geometry (in the Brill-Noether-theoretic sense)  by the measure of how far Castelnuovo and Mumford's (and Green's) theorems are from being optimal. 
Concretely, the Green-Lazarsfeld's normal generation theorem (\cite{gl-inv}) can be stated as follows: given a pair $(C,L)$  with $C$ of genus $g$ and $\deg L=2g+1-k$, where $k$ is a positive integer,  the module $\mathcal R_{C,L}(L)$ fails to be generated in degree $0$ if and only if either $L$ is not very ample or $\mathrm{cliff}(C)+2h^1(C,L)< k$.  Here $\mathrm{cliff}(C)$ denotes the Clifford index of the curve, which is a measure of how exceptional is the curve from the point of view of special divisors, see e.g. \cite{gl-inv}.

Concerning abelian varieties, we first remark that, not surprisingly, Koizumi's theorem (as well as its generalization to normal presentation, see \cite[Theorem 6.13]{kempf-libro}, and to syzygies, \cite{sav}, see also \cite{caucci}, \cite{caucci-kummer}) is usually not optimal for non-principal polarizations, and it needs to be improved in function of certain measures of the positivity of the polarization, which in turn depend on the geometry of the polarized abelian variety $(A,L)$.  In recent years there has been considerable research on this (\cite{iyer}, \cite{h-t}, \cite{lpp},\cite{jp}, \cite{kl},\cite{lo},\cite{caucci}, \cite{ito3}, \cite{ito4},\cite{ito-3-folds}, \cite{zhi-survey}, \cite{zhi-forum} and references therein). However this does not apply to what is arguably the most important case, namely  principally polarized abelian varieties $(A,\OO_A(\Theta))$, where in any case Koizumi's theorem cannot be improved to a better statement concerning the generation in degree $0$ of modules of the form (\ref{eq:koizumi}). This is  because in the multiplication maps of global sections
\begin{equation}\label{eq:mult2}
H^0(A,\OO_A(2\Theta))\otimes H^0(A,\OO_A(2\Theta)\otimes P_\alpha)\rightarrow H^0(\OO_A(4\Theta)\otimes P_\alpha)
\end{equation}
the dimensions of the source and  of the target are equal, namely $4^g$, where $g=\dim A$. Hence such maps fail to be surjective for $\alpha$ belonging to a divisor on $\Pic0 A$ (which can be computed, see \S3 below). Therefore there is no way of distinguishing any p.p.a.v. by the failure of generation of degree $0$ of modules of the form (\ref{eq:koizumi}).  It can be similarly shown that the same is true for properties as having a linear minimal free resolution at stage $p$ over the symmetric algebra. More conceptually, we recall that  on abelian varieties properties as normal generation (and, possibly to a lesser extent, linearity of the minimal resolution up to a certain stage) are governed by the \emph{base point freeness threshold} (an invariant introduced in \cite[\S8]{jp} and extensively studied e.g. in \cite{caucci}, \cite{ito4}, \cite{ito-3-folds}, \cite{zhi-survey}),  and   all p.p.a.v.'s have the same  base point freeness threshold,  namely $1$.

\noindent \emph{Hyperelliptic curves and Theorem \ref{th:first}. } With Theorem \ref{th:first} we propose a way to go beyond Koizumi's theorem on p.p.a.v.'s, 
in the spirit of the Green-Lazarsfeld normal generation theorem.  To see this more precisely, let us go back to curves for a moment. The first case of Green-Lazarsfeld's normal generation theorem, namely, in the above notation, $k=1$, tells that: if $\deg L=2g$, then the module $\mathcal R_{C,L}(L)$ is not generated in degree $0$ if and only if $L$ is not very ample or $\mathrm{cliff}(C)=0$, i.e., by Clifford's theorem, $C$ is hyperelliptic (this specific result was previously anticipated by Lange-Martens in \cite{lm}). This was completed by a theorem of Eisenbud-Koh-Stillman (\cite[Theorem 2]{eks}) which, in particular, states that: if $L$ and $M$ are line bundles on $C$ such that $\deg L=\deg M=2g$ and the graded module $\mathcal R_{C,L}(M)$ fails to be generated in degree $0$, then $L=M$ (in which case Green-Lazarsfeld's theorem applies).  Summarizing: \emph{a curve of genus $g\ge 2$ is  hyperelliptic if and only if, \emph{for each} line bundle $L$ of degree $2g=(2g-2)+2$ on $C$ there is another line bundle $M$ of degree $2g = (2g-2)+3-1$  such that  the module $\mathcal R_{C,L}(M)$ is not generated in degree zero}  (actually $L=M$). Theorem \ref{th:first} can be seen as the p.p.a.v.'s analog of this characterization of hyperelliptic curves.

Of course the question is how to go further. To date we do not have a precise conjecture, but we point out that a natural way to extend Theorem \ref{th:first} would be to consider the non-generation in degree $0$ of graded modules of the form 
\begin{equation}\label{eq:conj?}
\mathcal R_{A,\OO_A(2\Theta)\otimes P_\alpha}(\I_Z(3\Theta)\otimes P_\beta)
\end{equation}
 with $Z$ a suitable collection of points or  $0$-dimensional subscheme of suitable length. For example in Subsection \ref{rem:rem}  it is shown that on $g$-dimensional jacobians, or polarized products of g-dimensional jacobians with another p.p.a.v., there are collections $Z$ of $g$ points, satisfying suitable genericity hypothesis, such that  the modules (\ref{eq:conj?}) are not generated in degree $0$ for some choices of $\alpha,\beta\in\Pic0 A$. Furthermore, if the jacobian is hyperelliptic, a stronger condition holds.

%%%%%%%%%%%%%%%%%%%%%%%%%%
\subsection{Second characterization. } This makes use of higher gaussian maps (see  \S\ref{sect:notation}(\ref{it:gaussian}) below). They were introduced, in the context of deformation theory of linear series on curves, in \cite[Chapter 9]{acgh} and,  in a general setting, in \cite{wahl}.  Here we will assume that the abelian varieties in question are just polarized, i.e. equipped with a polarization, not necessarily principal (but elliptic curves will be always equipped with their natural, i.e. principal, polarization).  

\begin{theoremalpha}\label{th:second} Let $(A,\OO_A(\Theta))$ be a polarized abelian variety \emph{(not necessarily principally polarized)}. Let $p$ an integer $\ge 2$. The following are equivalent:
\begin{enumerate} \item $A$ is isomorphic, as polarized variety, to $E\times B$, where $\dim E=1$.
\item There exists $\alpha\in \Pic0 A$ such that the  $p$-th order gaussian map 
\begin{equation}\label{eq:g(alpha)}
\mathrm{Rel}^p_A(\OO_A(2(p+1)\Theta),\OO_A(2(p+1)\Theta)\otimes P_\alpha)\rightarrow H^0(A,S^p\Omega^1_A\otimes\OO_A(4(p+1)\Theta)\otimes P_\alpha)\end{equation}
is not surjective.
\item The  gaussian map
\begin{equation}\label{eq:g}
\mathrm{Rel}^p_A(\OO_A(2(p+1)\Theta),\OO_A(2(p+1)\Theta))\rightarrow H^0(A,S^p\Omega^1_A\otimes\OO_A(4(p+1)\Theta))
\end{equation}
is not surjective.
\end{enumerate}
\end{theoremalpha} 

We remark that,  for principal polarizations, the theorem is false for $p=0,1$ (see the next Subsection). The crucial case is $p=2$, asserting that $A\cong E\times B$ if and only if there exists some $\alpha\in\Pic0 A$ such that 
\begin{equation}\label{eq:g2(alpha)}
\mathrm{Rel}^2_A(\OO_A(6\Theta),\OO_A(6\Theta)\otimes P_\alpha)\rightarrow H^0(A, S^2\Omega^1_A\otimes \OO_A(12\Theta)\otimes P_\alpha)
\end{equation}
is not surjective, or the analogous condition (3) holds. As well as Corollary \ref{cor:B} above, in the case of principal polarizations, this gives, at least locally, equations for $\mathcal A_1\times \mathcal A_{g-1}$ inside $\mathcal A_g$.

Again, to put things into perspective, it is useful to compare this result with a theorem of  Bertram-Ein-Lazarsfeld of the 80's  on gaussian maps on curves (\cite[Theorem 1.7]{bel}):

\noindent (a) given  a smooth curve $C$ of genus $g\ge 1$, for all pairs of  line bundles $L$ and $M$ on $C$ such that
$\deg L, \deg M\ge  (g+1)(p+1)$ and $\deg L+\deg M> (p+2)(2g-2)+4(p+1)$, the $p$-th order gaussian map
\begin{equation}\label{eq:g(L,M)}
\mathrm{Rel}^p_C(L,M)\rightarrow H^0(C,\omega_C^{\otimes p}\otimes L\otimes M)
\end{equation}
 is surjective, \\
  %(2) if $C$ is non-hyperelliptic the same is true if  $\deg L+\deg M\ge 4(2g-2)+12$,\\
   (b) for each $d$ such that $ (g+1)(p+1)\le d\le (p+2)(2g-2)+4(p+1)- (g+1)(p+1)$ there exist pairs of line bundles  $L$ and $M$ on $C$ with $\deg L=d$ and $\deg L+\deg M=(p+2)(2g-2)+4(p+1)$ (particular case: $\deg L = \deg M = (p+1)(g-1)+2(p+1)$) 
such that the map (\ref{eq:g(L,M)}) is not surjective if and only if $g(C)=1$ or $C$ is hyperelliptic.

 Part (a) of the above theorem was (partially) extended to abelian varieties in \cite[Theorem C]{beppe-compo} (in turn re-proved and improved in \cite[Theorem 6.1]{nelson}), stating that, if  $n,m\ge 2(p+1)$ and $n+m>4(p+1)$ then the $k$-th order  gaussian map 
\[ \mathrm{Rel}^p_A(\OO_A(n\Theta),\OO_A(m\Theta)\otimes P_\alpha)\rightarrow H^0(A, S^{p}\Omega^1_A\otimes \OO_A((n+m)\Theta)\otimes P_\alpha)
\]
is surjective for all  $\alpha\in \Pic0 A$. 
So Theorem \ref{th:second} can be seen as the abelian varieties analog of part (b)  (holding for $p\ge 2$). 

We point out that there is also a more general formulation, allowing  in Theorem \ref{th:second}(b) or (c) two different multiples of the polarization, but one usually needs to replace line bundles with semihomogenous vector bundles, see Subsection \ref{sub:variants2}.

\subsection{p-jets separation}\label{sub:jets} Contrary to Theorem \ref{th:first}, the proof of Theorem \ref{th:second} is not direct (in fact it would be interesting to have a direct proof). In fact what we prove is  an equivalent characterization in terms of jets separation: 
\begin{theoremalpha}\label{th:third} Let $(A,\OO_A(\Theta))$ be a polarized abelian variety. The following are equivalent:
\begin{enumerate} 
\item $A$ is not isomorphic, as polarized variety, to $E\times B$, where $\dim E=1$.

\item There exists an integer $p\ge 2$ such that for all $\alpha\in\Pic0 A$ the line bundle $\OO_A((p+1)\Theta)\otimes P_\alpha$ is $p$-immersive, \emph{(namely separates $p$-jets at every point of $A$, i.e. the natural map
\[
H^0(A,\OO_A((p+1)\Theta)\otimes P_\alpha)\rightarrow H^0(A, \OO_A(p+1)\Theta)\otimes P_\alpha\otimes (\OO_A/\I_x^{p+1}))
\]
is surjective for all $x\in A$)}.

\item For all integers $p\ge 2$ and all $\alpha\in\Pic0 A$ the line bundle $\OO_A((p+1)\Theta)\otimes P_\alpha$ is $p$-immersive.

\item For all $\alpha\in\Pic0 A$ the line bundle $\OO_A(3\Theta)\otimes P_\alpha$ is $2$-immersive.
\end{enumerate}
\end{theoremalpha}
Note that, for principal polarizations, this is not true for $p=0,1$ because $\OO_A(\Theta)$ is not base point free and $\OO_A(2\Theta)$ is not immersive.

The proof of the equivalence between Theorem \ref{th:third}  and Theorem \ref{th:second} makes use of work of the first author (\cite{nelson}), recalled in full detail in Subsection \ref{sub:thresholds} below. Briefly, given a polarization $L$, one can make sense of the notion of $p$-immersivity of a \emph{fractional multiple} of $L$, as well as the notion of surjectivity, \emph{for all $\alpha\in\Pic0 A$}, of gaussian maps of order $p$ involving fractional powers of line bundles . This leads to notions of \emph{$p$-jets separation thresholds}, denoted $\epsilon_p(\l)$,  and of \emph{$p$-higher gaussian maps  surjectivity thresholds}, denoted $\mu_p(L)$. For every $p$ the Fourier-Mukai transform establishes a precise formula between the two thresholds, see (\ref{eq:thresholds}) below. The equivalence between theorems \ref{th:second} and \ref{th:third} follows from this. 

As for Theorem \ref{th:first}, the most interesting question is how to go further. In this case there is a clear (conjectural) answer, related to the Seshadri constant of the polarized abelian variety, denoted here $\epsilon_A(\utheta)$. In fact we expect that  for a $g$-dimensional polarized abelian variety, to satisfy upper bounds  
\begin{equation}\label{eq:bound}
\epsilon_A(\utheta)\le \lambda, 
\end{equation}
should be equivalent to satisfy conditions of non-surjectivity of certain gaussian maps of suitable higher order $p_0$ (depending on $g$ and $\lambda$),  or, equivalently, suitable  non-$p_0$-immersivity conditions for (possibly fractional) multiples of a polarization. More precisely, in the above framework of thresholds, to satisfy (\ref{eq:bound}) should be equivalent, for a certain positive integer $p_0$,  to  lower bounds on the $p_0$-jets separation thresholds 
\begin{equation}\label{eq:expectation}
 \epsilon_{p_0}((p_0+1)\l)\ge \frac 1 \lambda .
\end{equation}
This is motivated by the results in \cite[\S4]{nelson}, where  the first author developed, in the context of abelian varieties, a version of Demailly's theorem relating the Seshadri constant to asymptotic $p$-jets separation (\cite[Theorem 5.1.17]{pag1}, see also \cite{bs}),   in terms of the afore mentioned thresholds. Moreover, in \emph{loc. cit} the author made a finiteness conjecture whose veracity would give the existence of $p_{0}$ as in (\ref{eq:expectation}). Since it is known that  lower bounds as  (\ref{eq:expectation}) are equivalent  to suitable lower bounds on surjectivity threshold for gaussian maps of order $p_0$ (see Subsection \ref{sub:sesh}), this would produce some sort of explicit equations for the locus of polarized abelian varieties whose Seshadri constant is bounded by a given value.

All this is explained in Subsection \ref{sub:sesh} below, where it is also explained that Theorem \ref{th:third} can be seen as an effective version of Nakamaye's theorem characterizing polarized  abelian varieties  of the form $E\times B$ as the only ones achieving the minimal Seshadri constant (namely $\epsilon_A(\utheta)=1$, \cite[Theorem 1.1]{nakamaye}), and confirms the above mentioned finiteness conjecture of the first author. 

Let us exemplify all this in the interesting case of hyperelliptic jacobians. A well known conjecture of Debarre  (\cite{debarre-seshadri}) predicts that an \emph{irreducible} p.p.a.v. $A$ is such that  $\epsilon_A(\utheta)<2$ if and only if it is a hyperelliptic jacobian. Via the above process, this is translated in the following conjecture: there is a positive integer $p_0$  (possibly depending on $g:=\dim A$) such that  an i.p.p.a.v. $(A,\OO_A(\Theta))$ is a hyperelliptic jacobian 
if and only if 
 the $p_0$-th order gaussian map 
\[Rel^{p_0}(\OO_A((p_0+1)\Theta),\OO_A((p_0+1)\Theta)\otimes P_\alpha)\rightarrow H^0(S^{p_0}\Omega^1_A\otimes \OO_A(2(p_0+1)\Theta)\otimes P_\alpha)
\]
fails to be surjective  for all $\alpha\in\Pic0 A$.

  %%%%%%%%%%%%%%%%%%%%%%%%%%%%%%%
\bigskip
\noindent
  \section{ Notation, terminology and basic facts used in the paper}\label{sect:notation}

We will work over an algebraically closed ground field $k$ of characteristic zero.

\begin{enumerate} 

\item\emph{Kernel bundle. }  For a base point free line bundle $L$ on the abelian variety $A$ we denote $M_L$ the locally free sheaf defined by the exact sequence
 \[
 0\rightarrow M_L\rightarrow H^0(A,L)\otimes\OO_A\rightarrow L\rightarrow 0.
 \]

\item\emph{Polarizations, isogenies and translations. }\label{it:isog} Given an ample line bundle $L$ on an abelian variety $A$, we denote $\l$ its class in the N\'eron-Severi group and 
\[
\varphi_{\l}:A\rightarrow \Pic0 A
\]
the corresponding isogeny.  
Given a non zero integer $k$, 
\[
k_A:A\rightarrow A
\]
 denotes the multiplication-by-$k$ isogeny $x\mapsto kx$, and $t_x:A\rightarrow A$ the translation by a point $x\in A$. 

\medskip
 \item\emph{Fourier-Mukai functors. } A Fourier-Mukai functor between bounded derived categories, associated to a kernel $\F\in D(X\times Y)$ is denoted 
  \[
  \Phi_\mathcal{F}:D(X)\rightarrow D(Y).
  \]

 \item\emph{Cohomological support loci. } Given an an abelian variety $A$, we denote $\cP$ the Poincar\'e line bundle on $A\times \Pic0 A$. Given $\alpha\in\Pic0 A$ (respectively $x\in A$) the line bundle on $A$ (resp. on $\Pic0 A$) parametrized by a closed point $\alpha$ (resp. by $x$) is denoted $P_\alpha$ (resp. $P_x$). Given a coherent sheaf $\F$ on $A$, its \emph{cohomological support loci} are the subvarieties (of $\Pic0 A$)
 \[V^i(A,\F)=\{\alpha\in\Pic0 A\>|\> h^i(A,\F\otimes P_\alpha)>0\}.
 \]

 \item\emph{Vanishing conditions on coherent sheaves. }\label{it:van} 
  A coherent sheaf $\F$ on an abelian variety $A$ is said to \emph{satisfy the index theorem with index zero}, or, for short,  to be \emph{IT(0)}, if
 $V^i(A,\F)=\emptyset$ for all $i>0$. It is said to \emph{satisfy generic vanishing}, or, for short, to be a \emph{GV sheaf}, if 
 \[\mathrm{codim}_{\Pic0 A} V^i(A,\F)\ge i
 \]
  for all $i>0$. It is known by Hacon  and Pareschi-Popa (\cite{hacon}, \cite{pp-gv}) that $\F$ is a GV sheaf  if and only if its derived dual $\F^\vee:=R\mathcal Hom(\F,\OO_A)$ \emph{satisfies the weak index theorem with index $g$} (where $g=\dim A$), i.e. its Fourier-Mukai transform (with kernel $\cP^\vee$) is a sheaf in cohomological degree $g$: 
 \begin{equation}\label{eq:WIT}
 \Phi_{\cP^\vee}(\F^\vee)=R^g\Phi_{\cP^\vee}(\F^\vee)[-g].
 \end{equation}
  Finally, $\F$ is said to be \emph{M-regular} if 
  \[
  \mathrm{codim}_{\Pic0 A} V^i(A,\F)> i
  \] 
  for all $i>0$. This is equivalent to the fact that (\ref{eq:WIT}) holds and the sheaf $R^g\Phi_{\cP^\vee}(\F^\vee)$ is torsion free.
  Needless to say, a coherent sheaf $\F$ satisfies conditions IT(0), or GV or M-regularity if and only if $\F\otimes P_\alpha$ does, for every $\alpha\in\Pic0 A$. 
  
 \medskip \item\emph{$\mathbb Q$-twisted coherent sheaves. } According to \cite[Definitions 6.2.1 and 6.2.2]{pag2} $\mathbb Q$-twisted coherent sheaves  are equivalence classes of pairs $(\F,\udelta)$, where $\F$ is a coherent sheaf on an abelian variety $A$ and $\udelta$ is a class in $\mathrm{NS}(A)_{\mathbb Q}$, under the equivalence relation generated by $(\F\otimes L,\delta)\sim (\F,\l+\udelta)$ (where $L$ is an ample  line bundle and $\l$ denotes its N\'eron-Severi class). The equivalence class of $(\F,\udelta)$ is denoted $\F\langle\udelta\rangle$. 
  
\medskip \item\emph{Vanishing conditions on $\mathbb Q$-twisted coherent sheaves. }\label{it:Q-van}  According to \cite[\S5]{jp} the above vanishing conditions extend as follows:  keeping the notation above, given a rational number $\lambda=\frac a b$, the $\mathbb Q$-twisted coherent sheaf $\F\langle\lambda\l\rangle$ is said to be IT(0) (respectively GV, M-regular) if the coherent sheaf $b_A^*\F\otimes L^{\otimes ab}$ is (see (\ref{it:isog}) above for the notation). As a consequence of the main result of \cite{hacon} we have that:\\
(a)  to be IT(0) is an open condition, in the sense that \emph{$\F\langle\lambda\l\rangle$ is IT(0) if and only if for all sufficiently small  $\varepsilon \in\mathbb Q^+$,  $\F\langle(\lambda-\varepsilon)\l\rangle$ is IT(0)} (\cite[Theorem 5.2(c)]{jp}); \\
(b) on the other hand \emph{$\F\langle\lambda\l\rangle$ is GV if and only of $\F\langle(\lambda+\varepsilon)\l\rangle$ is IT(0) for all $\epsilon\in\mathbb Q^+$} (\cite[Theorem 5.2(a)]{jp}). 
 
\medskip \item \emph{Gaussian maps}\label{it:gaussian} (see e.g. \cite[\S9]{acgh}, \cite{wahl}, \cite[\S1]{bel},   \cite[\S1]{beppe-compo}). 
 On a smooth projective variety $A$, let $L$ and $M$ two line bundles and let $\Delta\subset A\times A$ be the diagonal.
 The $p$-th order gaussian map is the linear map
 \[\gamma^p_{A,L,M}:H^0(A\times A, (L\boxtimes M)\otimes \I_{\Delta}^p)\rightarrow 
  H^0(A\times A, ( (L\boxtimes M))\otimes \I_{\Delta}^p/\I_{\Delta}^{p+1})=H^0(A, L\otimes M\otimes S^p\Omega^1_A)
  \]
  It is customary to refer to the source as \emph{the space of $p$-th order relations between $L$ and $M$} and denote it
    \[
 \mathrm{Rel}^p_A(L,M):=H^0(A\times A,  (L\boxtimes M)\otimes\I_\Delta^p).
 \]
  \end{enumerate}

  %%%%%%%%%
 \section{Proof of Theorem \ref{th:first}}
 
 In the first place we note that Corollary \ref{cor:B} follows immediately from Theorem \ref{th:first} because on elliptic curves the map (\ref{eq:0}) is certainly not surjective, hence the same happens for decomposable p.p.a.v's of the form $E\times B$.  In turn, the proof of Theorem \ref{th:first} is made of some steps, as follows.

 \subsection{First reduction} 
 
 The main statement we will prove is the following 
 
  \begin{theorem}\label{th:reduction} If $\dim A\ge 2$ and $\Theta$ is irreducible then the vector bundles $M_{\OO_A(2\Theta)}\otimes
  \OO_A(n\Theta)\otimes P_\alpha $ are globally generated for all $n\ge 3$ and $\alpha\in\Pic0 A$. 
 \end{theorem}

 This implies Theorem \ref{th:first}  by the following  standard argument. In the first place we recall that $\Theta$ is reducible if and only if  the p.p.a.v. $(A,\OO_A(\Theta))$  is \emph{decomposable}, namely the product of  lower dimensional p.p.a.v.'s  $(A_i,\OO_{A_i}(\Theta_i))$ (e.g. \cite[Theorems 4.3.1-6]{birke-lange} or \cite[Corollary 10.4]{kempf-libro}). If this is the case then the graded module \[\mathcal R_{A,\OO_A(2\Theta)}(\I_x(3\Theta)\otimes P_\alpha)\] is generated in degree zero if and only if the same is true on each factor $A_{i}.$ Therefore (1) implies  (2) (and (3)), because (2) is obviously true on elliptic curves, since, given a point $p$ on an elliptic curve $E$, the multiplication map $H^0(E,\OO_E(2p))\otimes H^0(E,\OO_E(2p))\rightarrow H^0(E,\OO_E(4p))$ is obviously not surjective. It remains to prove that (3) implies (1). As above it is sufficient to prove that, under the assumption that $\dim A\ge 2$ and $\Theta$ irreducible, the multiplication maps 
 \begin{equation}\label{eq:mult} S^nH^0(A,\OO_A(2\Theta))\otimes H^0(A,\I_x(3\Theta)\otimes P_\alpha)\rightarrow
 H^0(A, \I_x((3+2n)\Theta)\otimes P_\alpha)
 \end{equation}
  are surjective for all $n\ge 1$, $x\in A$ and $\alpha\in\Pic0 A$. 

To prove the surjectivity of these maps we consider the following commutative diagram
\[
\xymatrix{H^0(\OO_A(2\Theta))\otimes S^{n}H^0(\OO_A(2\Theta))\otimes H^0(\I_x(3\Theta)\otimes P_{\alpha})\ar[r]\ar[d] & S^{n+1}H^{0}(\mathcal{O}_{A}(2\Theta))\otimes H^{0}(I_{x}\otimes\mathcal{O}_{A}(3\Theta)\otimes P_{\alpha}) \ar[d]\\
H^0(\OO_A(2\Theta))\otimes H^0(\I_x((3+2n)\Theta)\otimes P_\alpha)\ar[r]&  H^0(\I_x((3+2(n+1))\Theta)\otimes P_\alpha).\\
}
 \]
From this diagram we see, inductively, that in order to prove the surjectivity of the maps \eqref{eq:mult}, it suffices to prove the surjectivity, for all $n\ge 0,$ of the multiplication maps of global sections
 \begin{equation}\label{eq:mult2} H^0(A,\OO_A(2\Theta))\otimes H^0(A,\I_x(3+2n)\Theta)\otimes P_\alpha)\rightarrow H^0(A,\I_x(3+2(n+1))\Theta)\otimes P_\alpha).
 \end{equation}

 Now, to prove the surjectivity of \eqref{eq:mult2}, we start by noticing that the sheaves $\I_x((3+2n)\Theta))$  are IT(0) for all $n\ge 0,$ because the line bundles $\OO_A((3+2n)\Theta)\otimes P_\alpha$  are base point free for all $n\ge 0$ and $\alpha\in\Pic0 A$. Therefore, taking cohomology in the exact sequence
 \[0\rightarrow M_{\OO_A(2\Theta)}\otimes \I_x((3+2n)\Theta))\otimes P_\alpha\rightarrow H^0(\OO_A(2\Theta))\otimes \I_x((3+2n)\Theta))\otimes P_\alpha\rightarrow
 \I_x((3+2(n+1)\Theta)\otimes P_\alpha\rightarrow 0,
 \]
 we get that  the surjectivity of the multiplication maps of global sections (\ref{eq:mult2})  for all $x\in A$ and all $\alpha\in \Pic0A$ is equivalent to the fact that the coherent sheaves  $M_{\OO_A(2\Theta)}\otimes \I_x((3+2n)\Theta)$ are IT(0) for all $x\in A$. 
 The standard exact sequence
\[0\rightarrow  M_{\OO_A(2\Theta)}\otimes \I_x((3+2n)\Theta)\rightarrow M_{\OO_A(2\Theta)}\otimes \OO_A((3+2n)\Theta)\rightarrow (M_{\OO_A(2\Theta)}\otimes \OO_A((3+2n)\Theta))\otimes k(x)\rightarrow 0,
\]
where $k(x)\cong k$ is the residue field at $x,$ shows that the sheaves $M_{\OO_A(2\Theta)}\otimes \I_x((3+2n)\Theta)$ are IT(0) for all $x\in A$ as soon as $M_{\OO_A(2\Theta)}\otimes \OO_A((3+2n)\Theta)$ is IT(0) and globally generated. But it is well known that $M_{\OO_A(2\Theta)}\otimes \OO_A((3+2n)\Theta)$ is IT(0)  for $n\ge 0$ (this is  a particular case of Koizumi's theorem recalled in the introduction). This proves that Theorem \ref{th:reduction} implies the surjectivity of the maps (\ref{eq:mult2}) for all $n\ge 0$,  $x\in A$ and $\alpha\in\Pic0 A$ and hence it concludes the proof that Theorem \ref{th:first} follows from Theorem \ref{th:reduction}.

 \subsection{Second reduction} Turning to the proof of Theorem \ref{th:reduction}, the \emph{generation} properties (in the sense of the paper \cite{p4}) of the vector bundle $M_{\OO_A(2\Theta)}\otimes\OO_A(2\Theta)$ play an essential role. In fact the main point is
 \begin{lemma}\label{lem:generation} 
 Let $(A,\utheta)$ be an indecomposable p.p.a.v. of dimension $g\ge 2$. \begin{enumerate} 
 
 \item The subvariety $V^0(A,M_{\OO_A(2\Theta)}\otimes\OO_A(2\Theta))$ (with its reduced structure) is a prime  divisor of $\Pic0 A$ algebraically equivalent to $4^{g-1}\widehat{\utheta}$ \emph{ (where $\hat\utheta$ denotes the principal polarization on the dual abelian variety $\Pic0 A$)}.
 
 \item  the map
\[
 \mathrm{ev}_U:\bigoplus_{\beta\in U}H^0(A, M_{\OO_A(2\Theta)}\otimes\OO_A(2\Theta)\otimes P_\beta)\otimes P_\beta^{\vee}\rightarrow M_{\OO_A(2\Theta)}\otimes\OO_A(2\Theta)
\]
 is surjective for all nonempty open subsets $U\subset \Pic0 A$ meeting $V^0(A,M_{\OO_A(2\Theta)}\otimes\OO_A(2\Theta))$. 
 \end{enumerate}
 \end{lemma} 
 \noindent (we recall that, in the terminology introduced in (\cite[Definition 1.1]{p4}) statement (2) means that the vector bundle  $M_{\OO_A(2\Theta)}\otimes\OO_A(2\Theta)$ is \emph{generated by the set consisting of the single irreducible subvariety $\{V^0(A,M_{\OO_A(2\Theta)}\otimes\OO_A(2\Theta))\}$}).

 Let us show that Lemma \ref{lem:generation} implies Theorem \ref{th:reduction}. This follows at once from \cite[Proposition 2.3.1]{p4}. For sake of self-containedness, let us reproduce it in the case at hand. For $x\in A$ and $\alpha\in\Pic0 A,$ let $U^\alpha_x$ be the open set of $\Pic0 A$ defined as $U^\alpha_x=\{\beta\in\Pic0 A\>|\> x\not\in \Theta_{\alpha-\beta}\}$, where, for $\gamma\in\Pic0 A$, $\Theta_{\gamma}$ denotes the divisor corresponding to the line bundle $\OO_A(\Theta)\otimes P_\gamma $ (which is the translate of $\Theta$ by $\gamma$, where $\gamma$ is seen as a point of $A$ via the identification induced by the principal polarization). We have that $x\in \Theta_{\alpha-\beta}$ if and only if $\beta\in \Theta_{\alpha-x}$. Therefore $U^\alpha_x$ is the complement of a theta divisor of $\Pic0 A$ and therefore it follows from  (1) of the Lemma that $U^\alpha_x$ must meet the irreducible subvariety $V^0(A,M_{\OO_A(2\Theta)}\otimes\OO_A(2\Theta))$ as soon as $g\ge 2$. Now, by definition, for each $\beta\in U_{x}^{\alpha},$ we have that the restriction map $H^{0}(A,\mathcal{O}_{A}(\Theta_{\alpha-\beta})) \rightarrow \mathcal{O}_{A}(\Theta_{\alpha-\beta})\otimes k(x)$
is surjective. Therefore, tensoring with $H^{0}(A,M_{\OO_A(2\Theta)}\otimes\OO_A(2\Theta)\otimes P_\beta)$ and taking the direct sum over $\beta\in U_{x}^{\alpha}$ we get a surjective map  
  \[ 
 \xymatrix{\bigoplus_{\beta\in U^\alpha_x}H^0(A, M_{\OO_A(2\Theta)}\otimes\OO_A(2\Theta)\otimes P_\beta)\otimes H^0(A,\OO_A(\Theta_{\alpha-\beta}))\ar[d]\\  \bigoplus_{\beta\in U^\alpha_x}H^0(A, M_{\OO_A(2\Theta)}\otimes\OO_A(2\Theta)\otimes P_\beta)\otimes \OO_A(\Theta_{\alpha-\beta})\otimes k(x).}
 \]

Now, composing with $\mathrm{ev}_{U^\alpha_x}$ (twisted with $\OO_A(\Theta)\otimes P_{\alpha}\otimes k(x)$) we get a map
 \[ 
 \bigoplus_{\beta\in U^\alpha_x}H^0(A, M_{\OO_A(2\Theta)}\otimes\OO_A(2\Theta)\otimes P_\beta)\otimes H^0(A,\OO_A(\Theta)\otimes P_{\alpha-\beta})\rightarrow M_{\OO_A(2\Theta)}\otimes\OO_A(3\Theta)\otimes P_{\alpha}\otimes k(x) ,
 \] 
 which is surjective by (2) of the Lemma. It follows that the vector bundle $\mathcal M:=M_{\OO_A(2\Theta)}\otimes\OO_A(3\Theta)\otimes P_\alpha$ is globally generated at $x$, because the last map factors through the evaluation of global sections  of $\mathcal M$ at $x$. Since $x$ is arbitrary, $\mathcal M$ is globally generated. This proves that Lemma \ref{lem:generation} implies Theorem \ref{th:reduction}.

\subsection{Proof of Lemma \ref{lem:generation}}  (1) Applying the Fourier-Mukai transform $\Phi_{\cP}$ to the twisted evaluation map 
\[
e: H^0(A,\OO_A(2\Theta))\otimes \OO_A(2\Theta)\rightarrow \OO_A(4\Theta)
\]
 we get a map of vector bundles on $\Pic0 A$ (both of rank $4^g$)
 \[
 \Phi_{\cP}(e): H^0(A,\OO_A(2\Theta))\otimes \Phi_{\cP}(\OO_A(2\Theta))\rightarrow \Phi_{\cP}(\OO_A(4\Theta))
 \]
 which, by base-change,  is fiberwise  the multiplication-of-global-sections map
 \[m_\alpha: H^0(A,\OO_A(2\Theta))\otimes H^0(A,\OO_A(2\Theta)\otimes P_\alpha)\rightarrow H^0(A,\OO_A(4\Theta)\otimes P_\alpha).
 \]
 Therefore the subvariety
 \[
 V^0(A,M_{\OO_A(2\Theta)}\otimes\OO_A(2\Theta))
 \] 
 (with the reduced scheme structure)  coincides with the support of the divisor of zeroes  
 \begin{equation}\label{eq:D}
 D:=(\det \Phi_{\cP}(e))_0.
 \end{equation}
  We claim that the N\'eron-Severi class of $D$ is $4^{g-1}\widehat{\utheta}$. The claim follows from the following standard computation of the class of the source and the target of $\Phi_{\cP}(e)$. First, we observe that, since $h^{0}(A,\OO_A(2\Theta)) = 2^{g},$ we have that
\begin{equation}\label{eq:classD}
[D] = [\det\Phi_{\mathcal{P}}(\OO_A(4\Theta))] - 2^{g}[\det\Phi_{\cP}(\OO_A(2\Theta))] 
\end{equation}
where $[ - ]$ stands for the class inside the Ner\'on severi of the corresponding divisor or line bundle. By \cite[Proposition 3.11(1)]{mukai-FM}, for any positive integer $n,$ we have that 
\begin{equation}\label{eq:mukailb}
\varphi_{n\theta}^{\ast}\Phi_{\cP}(\mathcal{O}_{A}(n\Theta))\cong H^{0}(A,\OO_A(n\Theta))\otimes\mathcal{O}_{A}(-n\Theta)\cong\OO_A(-n\Theta)^{\oplus n^{g}}
\end{equation}
On the other hand, we have that $\varphi_{n\utheta} = [n]_{\Pic0 A}\circ\varphi_\utheta,$ where $[n]_{\Pic0 A}$ denotes the multiplication-by-$n$ in $\Pic0 A$. Since $[n]_{\Pic0 A}^{\ast}$ acts as multplication by $n^{2}$ on the N\'eron-Severi group, combinig with \eqref{eq:mukailb} we get that 
\[\varphi_{\utheta}^{\ast}[\det\Phi_{\cP}(\mathcal{O}_{A}(n\Theta))] =-n^{g-1}\utheta = -\varphi_{\utheta}^{\ast}(n^{g-1}\widehat{\utheta})\] 
and thus (since $\varphi_\utheta$ is an isomorphism), we conclude that 
\[ [\det\Phi_{\cP}(\mathcal{O}_{A}(n\Theta))] = -n^{g-1}\widehat{\utheta}.\]
Finally, substituting in \eqref{eq:classD} we obtain
\[ [D] = (-4^{g-1}+2^{2g-1})\widehat{\utheta} = 4^{g-1}\widehat{\utheta},\]
as claimed.
 
 Now we consider the pulled back map
  \[\varphi_{2\utheta}^*
 \Phi_{\cP}(e): H^0(\OO_A(2\Theta))\otimes \varphi_{2\utheta}^*\Phi_{\cP}(\OO_A(2\Theta))\rightarrow \varphi_{2\utheta}^*\Phi_{\cP}(\OO_A(4\Theta))
 \]
 As above, the class of the divisor of zeroes of the determinant of this map, say $E$, is $4^g\utheta$. 
 
 If $\utheta$ is irreducible, assuming, as we can, that the divisor $\Theta$ is symmetric,  we have that 
 \begin{equation}\label{eq:E}
 E=\sum_{a\in A[2]}\Theta_{a}
\end{equation}
  (note that, since $\Theta$ is assumed to be irreducible, this is a sum of distinct prime divisors). This is shown in   \cite[Subsection 3.4]{ps}. Briefly, the proof is as follows. Fiberwise the map  $\varphi_{2\utheta}^*\Phi_{\cP}(e) $
 is the multiplication map of global sections
 \begin{equation}\label{eq:mx}
 m_x: H^0(A,\OO_A(2\Theta))\otimes H^0(A,\OO_A((2\Theta)_{x}))\rightarrow H^0(A,\OO_A(2\Theta+(2\Theta)_{x})).
 \end{equation}
 A theorem of Kempf (\cite[Theorem 3]{kempf-mult}, see \cite[Theorem 2.1]{ps} for the correct statement) says that the dimension of the cokernel of the map $m_x$ is the number of points in $A[2]\cap\Theta_x$, where $A[2]$ denotes the group of 2-division points of $A$. In particular the linear map $m_x$ is singular if and only if $x\in \Theta_a$ for some $a\in A[2]$ (recall that $\Theta$ is assumed to be symmetric).  
 This proves that the support of $E$ is the right hand side of (\ref{eq:E}). Since the class of the right hand side is equal to class of $E$, namely $4^g\utheta$, (\ref{eq:E}) is proved.
 
Finally, since $\varphi_{2\utheta}^*(\varphi_{2\utheta *}(\Theta))=\sum_{a\in A[2]}\Theta_{a}$, it follows that, up to translation by a $2$-division point, $D=\varphi_{2\utheta *}(\Theta)$. We recall that it is well known (see e.g. \cite[p.1591]{p3}) that $\varphi_{2\utheta *}(\Theta)$ is a prime divisor (of class $4^{g-1}$). Hence it coincides with the subvariety $V^0(A,M_{\OO_A(2\Theta)}\otimes\OO_A(2\Theta))$.

 \noindent (2) Let us denote, for typographical brevity, 
 \[
 \E:= M_{\OO_A(2\Theta)}\otimes\OO_A(2\Theta).
 \]
We note that $\E$ is a GV sheaf. To justify this, we start by noticing that by definition we have an exact sequence 
\begin{equation}\label{eq:defML}
0 \rightarrow \E\otimes P_{\alpha} \rightarrow H^{0}(\OO_A(2\Theta))\otimes \OO_A(2\Theta)\otimes P_{\alpha}\rightarrow\OO_A(4\Theta)\otimes P_{\alpha}\rightarrow 0
\end{equation}
and hence, since $H^{i}(\OO_A(2\Theta)\otimes P_{\alpha}) = H^{i}(\OO_A(4\Theta)\otimes P_{\alpha}) = 0$ for all $i\geq 1$ and $\alpha\in\Pic0 A,$ we see that $V^{i}(\E)$ is empty for all $i\geq 2.$ Moreover, we see that $\alpha$ lies in $V^{1}(\E)$ if and only if the multiplication map
\[H^{0}(\OO_A(2\Theta))\otimes H^{0}(\OO_A(2\Theta)\otimes P_{\alpha})\rightarrow H^{0}(\OO_A(4\Theta)\otimes P_{\alpha})\]
 is not surjective, while $\alpha$ lies in $V^{0}(\E)$) if and only if this map is not injective. Since the source and target of this map have the same dimension (namely, $4^{g}$), its surjectivity and injectivity are equivalent and thus $V^{1}(\E) = V^{0}(E).$ In particular, $V^{1}(\E)$ is a divisor and thus we conclude that $\E$ is a GV sheaf, as claimed. 

Now, let 
 \[\T:=\Phi_{\cP^\vee}(\E^\vee)[g]=R^g\Phi_{\cP^\vee}(\E^\vee),
 \]
where the equality holds because $\E$ is a GV sheaf, see \S\ref{sect:notation}(\ref{it:van}). Note that, by duality (\cite[Lemma 2.2]{pp-gv}), $\T$ is the (derived) dual of $\Phi_{\cP}(\E)$. 
 
 The key point is a general result of the second author (\cite[Corollary 6.2.1]{p4}) which, in the case at hand, states that the generation property for $\E$, i.e. the surjectivity of the evaluation maps $\mathrm{ev}_U$ as in Lemma \ref{lem:generation}(2), holds as soon as the following conditions hold: (a) the scheme-theoretic support of $\T$  is irreducible and reduced, and (b)  $\T$ is of pure dimension $g-1$. For sake of self-containedness we have included a simplified proof of this result in the case at hand in the subsection below.  
 
 Assuming this,  applying the functor $\Phi_{\cP^{\vee}}((-)^{\vee})$ to the exact sequence \eqref{eq:defML}, we note that, since $\E$ is GV, $\T$ sits in the locally free resolution 
 \[0\rightarrow R^g\Phi_{\cP^\vee}(\OO_A(-4\Theta))\rightarrow H^g(A,\OO_A(-2\Theta))\otimes R^g\Phi_{\cP^\vee}(\OO_A(-2\Theta))\rightarrow \T\rightarrow 0
 \]
 where the first map is the dual of $\Phi_{\cP}(e)$. Therefore $\T$ is of pure dimension $g-1$ (e.g. by \cite[Proposition 1.1.10(2)]{HL}). Moreover the scheme-theoretic support of $\T$ is the divisor $D$ of (\ref{eq:D}) which is irreducible and reduced. This proves  (2).

\subsection{Appendix: proof of the generation property} For the reader's convenience, we outline the proof of \cite[Corollary 6.2.1]{p4}  in the case at hand, which is much simpler than the general case. We need to show that the following hypotheses: \emph{(a)  the scheme-theoretic support of $\T$, namely the divisor $D$, is reduced and irreducible}, and  \emph{(b) $\T=R^g\Phi_{\cP^\vee}(\E^\vee)$ is pure of codimension $g-1$},  imply that  the sheaf $\E$ is generated, i.e. the evaluation maps $\mathrm{ev}_U$ as in Lemma \ref{lem:generation}(2) are surjective. 
  This  is equivalent to the condition that that the multiplication maps
 \begin{equation}\label{eq:ev(x)}
 \bigoplus_{\beta\in U}H^0(A, \E\otimes P_\beta)\otimes H^0(A, P_\beta^\vee\otimes k(x))\rightarrow H^0(A, \E\otimes k(x))
 \end{equation}
 (where $k(x)\cong k$ denotes the residue field at $x$) are surjective for all $x\in A$ and for all open sets $U\subset \Pic0 A$ meeting $D$.  
 
 Dualizing the individual maps of (\ref{eq:ev(x)}) we get 
 \[
 \mathrm{Ext}^g(k(x),\E^\vee)\rightarrow H^g(A, \E^\vee\otimes P_\beta^\vee)\otimes H^0(A,P_\beta^\vee\otimes k(x))^\vee .
 \]
 The Fourier-Mukai equivalence $\Phi_{\cP^\vee}:D(A)\rightarrow D(\Pic0 A)$ identifies the source of the above map to \[
 \mathrm{Hom}(P_x^\vee, \Phi_{\cP^\vee}(\E^\vee)[g])=\mathrm{Hom}(P_x^\vee, R^g\Phi_{\cP^\vee}(\E^\vee))
 \]
  (the equality follows from (\ref{eq:WIT}),  because, as noted above, $\E$ is a GV sheaf). 
 Moreover, by base-change, $H^g(A, \E^\vee\otimes P_\beta^\vee)$ is identified to fibre at $\beta$ of the coherent sheaf $R^g\Phi_{\cP^\vee}(\E^\vee)$ and 
 $H^0(A,P_\beta^\vee\otimes k(x))^\vee$ is identified to the fibre at $\beta$ of the line bundle $P_x$. In conclusion, the Fourier-Mukai equivalence
 $\Phi_{\cP^\vee}$ identifies the dual of the map (\ref{eq:ev(x)}) to the evaluation map
 \begin{equation}\label{eq:ev-dual} H^0(\Pic0 A, R^g\Phi_{\cP^\vee}(\E^\vee)\otimes P_x)\rightarrow 
 \prod_{\beta\in U} R^g\Phi_{\cP^\vee}(\E^\vee)\otimes P_x\otimes k(\beta)
 \end{equation}
 An element of the kernel of this map would be a global section vanishing identically on a non empty open subset of the support of the sheaf $\mathcal T=R^g\Phi_{\cP^\vee}(\E^\vee)$, which is pure of codimension 1, hence, seen as a sheaf on its (scheme theoretic) support, is torsion free. Moreover its scheme theoretic support  is reduced by hypothesis. Therefore such a global section must be zero.
 
 %%%%%%%%%%%%%%%%%%%%%%%%%%%%
 
 \section{2-jets separation and proof of Theorems \ref{th:second} and \ref{th:third}}
 
 \subsection{Proof of Theorem \ref{th:third}} We recall that a line bundle $L$ is said to separate $p$-jets at a point $x\in A$ if the map $H^0(A,L)\rightarrow H^0(A,L\otimes \OO_A/\I_x^{p+1})$ is surjective, and $p$-immersive if it separates $p$-jets at every point. We will prove the implications (1) $\Rightarrow$ (4) $\Rightarrow$ (3) $\Rightarrow$ (2) $\Rightarrow$ (1).

   In what follows, we will use repeatedly the fact that for any ample line bundle $L$ on an abelian variety $A$, we have that the set of translates $t_{x}^{\ast}L,$ $x\in A,$  coincides with the set of twists $L\otimes P_{\alpha},$ for $\alpha\in\Pic0 A.$ Moreover, two ample line bundles $L_{1},L_{2}$ have the same class in the N\'eron-Severi group of $A$ if and only if $L_{1}\cong L_{2}\otimes P_{\alpha}$ for some $\alpha\in\Pic0 A.$

   We start with a preliminary fact:

 \begin{lemma}\label{lem:M-regular}  The sheaf $\I^2_x(2\Theta)$ is M-regular \emph{(see \S\ref{sect:notation}(\ref{it:van}))} for all $x\in A$ and  if and only if $A$ is not isomorphic,  as polarized variety, to $ E\times B$, where $E$ is an elliptic curve.
 \end{lemma} 
 \proof  Since, as it is well known,  any ample class in the N\'eron-Severi group has a symmetric representative (e.g. \cite[Lemma 4.6.2]{birke-lange})\footnote{i.e a representative $L$ with  $(-1_{A})^{\ast}L\cong L,$ where $-1_{A}$ is the inversion map on $A$},
  we can assume that $\Theta$ is symmetric. Since $V^i(A,\I_x(2\Theta))$ is empty for $i\ge 2$, we have that $\I_x(2\Theta)$ is M-regular if and only if 
\[
\mathrm{codim}_{\Pic0 A} V^1(A,\I^2_x(2\Theta))>1.
\]
 Since all line bundles $\OO_A(2\Theta)\otimes P_\alpha$ are translates of $\OO_A(2\Theta)$ (and conversely) and $\I^2_x(t_y^*\OO_A(2\Theta))\cong t_y^*(\I^2_{x-y}(2\Theta)),$ we have that the locus
 $V^1(A,\I^2_x(2\Theta))$ has the same dimension as the locus $R=\{z\in A : h^1(A,\I^2_z(2\Theta))>0\}.$ Now, $R$ is nothing else than the ramification locus of the morphism $f:A\rightarrow \mathbb P(H^0(A,\OO_A(2\Theta)^\vee)$. Indeed: since $\OO_A(2\Theta)$ is globally generated and $H^{1}(A,\OO_A(2\Theta)) = 0,$ the condition $h^{1}(A,\I^{2}_z(2\Theta))\neq 0$ is equivalent to the non-separation of tangent vectors at $z$ and hence to the ramification of the associated morphism to the projective space. Moreover, if the polarized abelian variety  $(A,\OO_A(\Theta))$ is indecomposable then $R$ is either $0$-dimensional or empty. Indeed: by \cite[Theorem 4.5.5]{birke-lange} we have that $\OO_A(2\Theta)$ is very ample (in which case $R$ is empty) unless $\OO_A(\Theta)$ has divisorial base locus. Now, if $\OO_A(\Theta)$ has divisorial base locus and the polarized variety is indecomposable then the decomposition theorem (\cite[Theorem 4.3.1]{birke-lange}) ensures that $\OO_A(\Theta)$ defines a principal polarization, in which case, as it is well known, $R$ is just the set of $2$-torsion points of $A$ (e.g. \cite[Theorem 4.8.1]{birke-lange}) and hence it is zero-dimensional. Finally, if $(A,\OO_A(\Theta))$ is decomposed as the product of $(A_1,\OO_{A_1}(\Theta_1)\times (A_2,\OO_{A_2}(\Theta_2)),$ then the components of the ramification locus (if any) have dimension   equal to $\dim A_1$ or $\dim A_2$. In conclusion, the ramification locus has codimension $1$ if and only if one of the factors is $1$-dimensional. 
 \endproof
 
 Now we turn to the proof of Theorem \ref{th:third}.
 
 \noindent (1) $\Rightarrow$ (4). It is known that if  $\I^2_x(2\Theta)$ is M-regular, then $\I_x^2(3\Theta)$ is globally generated (this is an application of the basic M-regularity criterion of \cite[Theorem 2.4]{pp1}). In particular $\I_x^2(3\Theta)$ is globally generated at the point $x$ itself, i.e. the map $H^0(A,\I^2_x(3\Theta))\rightarrow H^0(A,(\I_x^2/\I_x^3)(3\Theta))$ is surjective, and it is easily seen that this means that $H^0(A,\OO_A(3\Theta))$ surjects onto 
 $H^0(A,(\OO_A/\I_x^3)(3\Theta))$, i.e. $\OO_A(3\Theta)$ separates $2$-jets at $x$. Since $x$ is arbitrary this proves the desired implication. 
 
\noindent   (4) $\Rightarrow $ (3) follows from the same argument (actually an easier version of it, because it follows that for $p\ge 3$ the sheaf  $\I_x^p(p\Theta)$ is IT(0), and not only M-regular).  

\noindent (2) $\Rightarrow$ (1). If $A\cong E\times B$ then no line bundle of the form $\OO_A((p+1)\Theta)\otimes P_\alpha)$ is $p$-immersive. Indeed, if $x=(x_{1},x_{2})\in E\times B$ then we have a surjective morphism $L\otimes\OO_A/I_{x}^{p+1}\rightarrow L\otimes\OO_{E\times x_2}/I_{x/E\times x_2}^{p+1}$ (of sheaves on $A$), where $I_{x/E\times x_2}$ is the ideal of $x$ inside $E\times x_2.$ The kernel of this morphism has $0$-dimensional support and hence this morphism remains surjective after taking global sections. On the other hand, we have a commutative diagram:
\[
\xymatrix{ H^{0}(A,L) \ar[r]\ar[d] & H^{0}\left(E\times x_2, \left.L\right|_{E\times x_2}\right)\ar[d] \\
               H^{0}\left(A,L\otimes\OO_A/I_{x}^{p+1}\right) \ar[r] & H^{0}\left(E\times x_2, L\otimes\OO_{E\times x_2}/I_{x/E\times x_2}^{p+1}\right).}\] 
In particular, if $L$ separates $p$-jets at $x$ we would have that $\left.L\right|_{E\times x_2}$ also separates $p$-jets at $x.$ However, $\left.\OO_A((p+1)\Theta)\right|_{E\times x_2}$ has degree $p+1$ and hence it can not separate $p$-jets at every point of the elliptic curve $E\times x_2.$ This concludes the proof of Theorem \ref{th:third}.

\subsection {p-jets separation thresholds}\label{sub:p-jets} In order to show the relation between Theorem \ref{th:third} and Theorem \ref{th:second}, in this subsection we recall the notions of \emph{$p$-jets separation thresholds} and \emph{$p$-th gaussian maps
surjectivity thresholds} introduced in the paper \cite{nelson}. 

Let $L$ be an ample line bundle on $A$. Since the set of all translates $t_y^*L$ coincides with the set of all line bundles of the form $L\otimes P_\alpha$, we have that $L$  is $p$-immersive if and only if, fixing a point $e\in A$ (for example the neutral element), the line bundles  $L\otimes P_\alpha$ separate $p$-jets at $e$ for all $\alpha\in\Pic0 A$. In turn, this is equivalent to the fact that the coherent sheaf $\I_e^{p+1}\otimes L$ is an IT(0) sheaf, as it follows from the fact that $H^{1}(A,L\otimes P_{\alpha}) = 0$ (for example, by Kodaira vanishing) and the long exact sequence in cohomology constructed from the definitional sequence 
\[ 0\to I_{e}^{p+1}\otimes L\otimes P_{\alpha}\to L\otimes P_{\alpha}\to L\otimes P_{\alpha}\otimes\OO_A/\I_{e}^{p+1}\to 0.\]  
Now, the IT(0) property can be defined also for $\mathbb Q$-twisted sheaves, as $\I_e^{p+1}\langle\lambda \l\rangle$ for $\lambda\in\mathbb Q$ (see \S\ref{sect:notation}(\ref{it:Q-van})). Therefore it is natural to define the \emph{$p$-jets separation thresholds} as
\[
\epsilon_p(\l)=\inf\{\lambda\in \mathbb Q\>|\> \I_e^{p+1}\langle\lambda \l\rangle \hbox{ is IT(0)}\}.
\]
 This generalizes to higher jets the \emph{base point freeness threshold} introduced in
 \cite[\S8]{jp} and \cite{caucci}, namely $\epsilon_0(\l)$. 
 
 By definition the $p$-jets thresholds are multiplicative, in the sense that, given a polarization $\utheta$ on $A$,  $\epsilon_p(k\utheta)=\frac 1 k \epsilon_p(\utheta)$. Moreover it is known (\cite[Corollary 3.9]{nelson}) that  
 \begin{equation}\label{eq:normalized}
 \epsilon_p((p+1)\utheta)=\frac{\epsilon_p(\utheta)}{p+1}\le 1
 \end{equation}
and   $\epsilon_p((p+1)\utheta)=1$, if and only if $(p+1)\utheta$ is not $p$-immersive. Hence Theorem \ref{th:third} can be rephrased as follows: \emph{let $p$ any fixed integer $\ge 2$. Then  $A$ is isomorphic, as polarized variety, to $E\times B$ if and only if $\epsilon_p((p+1)\utheta)=1$}.

 \subsection{Gaussian maps thresholds} Following a similar principle, in \cite[\S5]{nelson}  are defined \emph{higher gaussian maps surjectivity thresholds}. To recall this, we use the notation of \S\ref{sect:notation}(\ref{it:gaussian}) on gaussian maps. In the first place, applying ${p_2}_*$ to the exact sequence
 \[
 0\rightarrow \I_{\Delta}^{p+1}\otimes (L\boxtimes M)\rightarrow \I_{\Delta}^p\otimes(L\boxtimes M)\rightarrow (\I_{\Delta}^p/\I_{\Delta}^{p+1})(L\boxtimes M)\rightarrow 0
 \]
 one gets
 exact sequences
 \[0\rightarrow R^{p}_{A,L}\otimes M\rightarrow R^{p-1}_{A,L}\otimes M\rightarrow S^p\Omega^1_A\otimes L\otimes M
 \]
 where $R^i_{A,L}:={p_2}_*( \I_{\Delta}^{i+1}\otimes p_{1}^{\ast}L)$ and the gaussian map $\gamma^p_{L,M}$ is the $H^0$ of the last map. It is easily seen  that the short sequence above is exact also on the right as soon as the line bundle $L$ is $p$-immersive (see e.g. \cite[(5.3)]{nelson}), i.e., in the above notation, $\epsilon_p(L)<1$ (moreover if this is the case then the coherent sheaves $R^k_{A,L}$ are locally free for $k\le p$). Assuming that $\epsilon_p(L)<1$, the condition that the locally free sheaf $R^p_{A,L}$ is IT(0) implies that  the gaussian map $\gamma_{L,M\otimes P_\alpha}^p$ is surjective  for all $\alpha\in \Pic0 A$. Since the IT(0) condition makes sense for $\mathbb Q$-twisted sheaves, it is natural to define
 \[
 \mu_p(L)=\mathrm{inf}\{t\in\mathbb Q\>|\> R^p_{A,L}\langle tL\rangle \hbox{ is IT(0)}\}
 \footnote{The definition given in \cite[Definition 5.4]{nelson} is slightly different, but equivalent (see \cite[Lemma 5.3]{nelson}).}
 \]
  We have that $\mu_p(L)$ depends only on the N\'eron-Severi class of $L$, so one can call it $\mu_p(\l)$.  The fundamental relation with the $p$-jets separation threshold is   \cite[Theorem 5.5]{nelson}  asserting that, as soon as $\epsilon_p(\l)<1$, 
 \begin{equation}\label{eq:thresholds}
 \mu_p(\l)=\frac{\epsilon_p(\l)}{1-\epsilon_p(\l)}.
 \end{equation}
 
\begin{remark}\label{rem:gen} This is a generalization to arbitrary $p$ of the case $p=0$,  shown in \cite[Theorem D]{jp}, where the vector bundle $R^0_{A,L}$ is nothing else than the kernel bundle $M_L$ (\S\ref{sect:notation}(1)). In this case  the formula (\ref{eq:thresholds}), relates the surjectivity of multiplication maps of global sections threshold and the base point freeness threshold. This will be used in Subsection \ref{sub:variants1} below. 
\end{remark}

\subsection {Proof of Theorem \ref{th:second}}\label{sub:thresholds}  Using  (\ref{eq:thresholds}) Theorem \ref{th:second} follows  from (in fact is equivalent to)  Theorem \ref{th:third}. Indeed, as we saw, Theorem \ref{th:third} can be stated as follows: given an integer $p\ge 2$, $A$ is not isomorphic, as polarized variety,  to $E\times B$ if and only if $\epsilon_p((p+1)\utheta)<1$ i.e.  $\epsilon_p(2(p+1)\utheta)<\frac 1 2$. By (\ref{eq:thresholds}) this yields $\mu_p(2(p+1)\utheta)<1$ i.e. the $p$-th gaussian map 
\[
\label{eq:g(alpha)}
\mathrm{Rel}^p_A(\OO_A(2(p+1)\Theta),\OO_A(2(p+1)\Theta)\otimes P_\alpha)\rightarrow H^0(A,S^p\Omega^1_A\otimes\OO_A(4(p+1)\Theta)\otimes P_\alpha)
\]
 is surjective for all $\alpha\in\Pic0 A$. This proves that (2) implies (1) in Theorem \ref{th:second}. Now, it is enough to prove the implication (1) $\Rightarrow$ (3) in the case of elliptic curves. In fact the $p$-th  gaussian map (\ref{eq:g}) is never surjective on elliptic curves. This follows e.g. from \cite[Theorem C]{beppe-compo} , quoted in the introduction, or from and can be proved as follows.  Let $E$ be an elliptic curve and $L$ a line bundle of degree $2(p+1)$ on $E$. Then the gaussian maps
\begin{equation}\label{eq:elliptic}
\gamma^i_{E,L,L}: \mathrm{Rel}^i_E(L,L)\rightarrow H^0(E, S^p\Omega^1_E\otimes  L^{\otimes 2})=H^0(E, L^{\otimes 2})
\end{equation}
 are surjective for $i<p$. This is an easy particular case of   both \cite[Theorem C]{beppe-compo} and \cite[Theorem 1.7]{bel}, quoted in the introduction. Thus an easy 
calculation shows that $\dim \mathrm{Rel}^p_E(L,L)=4(p+1)$, i.e. the dimension of the source and of the target of the map (\ref{eq:elliptic})  are equal for $i=p$. But such a map cannot be injective because, as it follows from a well known computation in local coordinates (see e.g. \cite[\S 1.1 p.11-12]{lacopo} or \cite[\S 2]{ffl}, \cite[\S 5.1]{faro}), gaussian maps as (\ref{eq:elliptic}), i.e. $\gamma^i_{X,L,M}$ with $L=M$,  are zero on skew-symmetric (resp. symmetric) tensors) for $i$ even (resp. $p$ odd). This means the following. We have the filtration 
\[
\dots\subset \mathrm{Rel}^p(L,L)\subset \dots \subset\mathrm{Rel}^2(L,L)\subset \mathrm{Rel}^1(L,L)\subset H^0(L)\otimes H^0(L)
\]
(here we are neglecting the ambient variety in the notation). Then: \\
- $\mathrm{Rel}^1(L,L)=\wedge^2H^0(L)\otimes I^2(L)$, where $I^2(L)=\ker (S^2H^0(L)\rightarrow H^0(L^{\otimes 2}))$;  \\
- $\mathrm{Rel}^2(L,L)=\ker (\gamma^1_{L,L})_{|\wedge^2H^0(L)}\oplus I^2(L)$,\\
- $\mathrm{Rel}^3(L,L)=\ker (\gamma^1_{L,L})_{|\wedge^2H^0(L)}\oplus \ker(\gamma^2_{L,L})_{|I^2(L)}$,
 and so on. From this and the above computation it follows that  $\gamma^p_{E,L,L}$ is not injective, hence not surjective,  for $\deg L=2(p+1)$.

%%%%%%%%%%%%%%%%%%%%%%%%%%%%%

 \section{Remarks, variants and conjectures}\label{sect:remarks} 
 
 %%%%%%%%%%%%%%
  \subsection{A remark about Theorem \ref{th:first}} \label{rem:rem}
In view of possible extensions of Theorem \ref{th:first}, the conditions of being a jacobian, especially a hyperelliptic jacobian, or to be decomposed as the polarized product of those with an another p.p.a.v., seem to be relevant, as shown by the following generalization of the easy direction of Theorem \ref{th:first}.
 \begin{proposition} Let $A$ be the polarized product $J(C)\times B$, where $J(C)$ is the jacobian of a curve of genus $g$ and $B$ a p.p.a.v. . Let $Z\subset A$ be any collection of $g$ points $\{(z_1,b_1),\dots,(z_g,b_g)\}$ such that $\widetilde Z:=\{z_1,\dots,z_g\}$ is contained in (a translate of) $C$. Then
 \ there exist $\alpha,\beta\in\Pic0 A$ such that the $\mathcal S_{\OO_{J(C)}(2\Theta)\otimes P_\alpha}$-module 
  \[
  \mathcal R_{A,\OO_{A}(2\Theta)\otimes P_\alpha}(\I_{{Z}}(3\Theta)\otimes P_\beta)
  \]
   is not generated in degree $0$.
    If $C$ is hyperelliptic, for all $\alpha\in \Pic0 A$ there exists $\beta\in\Pic0 A$ such that the above holds. 
  \end{proposition} 
  \begin{proof} As usual, it is sufficient to assume that $A=J(C)$ and $Z$ a collection of $g$ points contained in $C$, and to show that there  exist $\alpha,\beta\in\Pic0 A$ such that the multiplication map 
 \begin{equation}\label{eq:jacobians}
 m_{\alpha,\beta}:H^0(A,\OO_{J(C)}(2\Theta)\otimes P_\alpha)\otimes H^0(A,\I_Z(3\Theta)\otimes P_\beta)\rightarrow H^0(A,\I_Z(5\Theta)\otimes P_{\alpha+\beta})
 \end{equation}
  is not surjective.  This map fits into the commutative diagram
   \begin{equation}\label{eq:jacobians2}
 \xymatrix{H^0(A,\OO_{J(C)}(2\Theta)\otimes P_\alpha)\otimes H^0(A,\I_Z(3\Theta)\otimes P_\beta)\ar[r]^{\>\>\>\>\>\>\>\>\>\>\>\>\>\>\>\>\>\>\>\>\>\>\>\>\>\>\>\>\>\>m_{\alpha,\beta}}\ar[d]& H^0(A,\I_Z(5\Theta)\otimes P_{\alpha+\beta})\ar[d]\\
 H^0(C,\OO_{C}(2\Theta)\otimes P_\alpha)\otimes H^0(C,O_{C}(3\Theta-Z))\otimes P_\beta)\ar[r]& H^0(C,\OO_{C}(5\Theta-Z))\otimes P_{\alpha+\beta})}
 \end{equation}
 It is clear that the vertical restriction map on the right is surjective. This is because its cokernel is contained in $H^1(A, \I_C(5\Theta)\otimes P_{\alpha+\beta}),$ which is zero because it is well known that $\I_{C}(5\Theta)$ is IT(0) (in fact much more is true, see \cite[Theorem 4.1]{pp1}).  
    Hence the non-surjectivity of the map $m_{\alpha,\beta}$ is implied by the non-surjectivity of the bottom  map of (\ref{eq:jacobians2}).  Choosing $\beta\in\Pic0 A$ ($\cong\Pic0 C$) such that 
    \begin{equation}\label{eq:eq}
    \OO_{C}(3\Theta-Z)\otimes P_\beta\cong \OO_{C}(2\Theta)\otimes P_\alpha
    \end{equation} the bottom map of (\ref{eq:jacobians2}) is
    \begin{equation}\label{eq:eq2}
    H^0(C,\OO_C(2\Theta)\otimes P_\alpha)\otimes H^0(C,\OO_C(2\Theta)\otimes P_\alpha)\rightarrow H^0(C,\OO_C(4\Theta)\otimes P_\alpha^{\otimes 2}).
    \end{equation}
    For $\alpha\in \Pic0 A$  such that the line bundle $\OO_{C}(2\Theta)\otimes P_\alpha$ is not very ample (the  line bundles of degree $2g$ on $C$ which are not very ample form a 2-dimensional family, namely $\omega_C(x+y)$, with $x,y\in C$),  
     one has the non-surjectivity of the map (\ref{eq:eq2}) (this is well known, for example it is part of the case $p=0$  of \cite[Theorem 2]{gl-compo}). Hence of the map $m_{\alpha,\beta}$ is not surjective. 

Now assume that the curve $C$ is hyperelliptic.   In this case for all $\alpha\in\Pic0 A$ \ the map (\ref{eq:eq2}) is not surjective (this is contained in the case $k=1$ of the general Green-Lazarsfeld normal generation theorem  recalled in the introduction, because the Clifford index of a smooth curve vanishes if and only if the curve is hyperelliptic). Therefore, for all $(\alpha,\beta)\in\Pic0 A\times\Pic0 A$  such that (\ref{eq:eq}) is satisfied the map $m_{\alpha,\beta}$ is not surjective.
   \end{proof}

  Finally, we remark  that  collections of points $Z$ in $J(C)$ such that $Z$ is contained in (a translate of) $C$ satisfy suitable genericity hypotheses as soon as $h^0(\OO_{C}(Z))=1$ (\cite[Example 3.4, Remark 6.5]{pp-cast}).

 \subsection{A variant of Theorem \ref{th:first}}\label{sub:variants1} One might think that Theorem \ref{th:first} depends on the  very  specific and interesting geometry of the maps   $A\rightarrow \mathbb P^{2^g-1}$ associated to the line bundles $\OO_A(2\Theta)$ (if the divisor $\Theta$ is symmetric this is the Kummer map). But in fact this is not the case, as it is shown by the following variants of  Theorem \ref{th:first}. These results do not seem to have a (known) analog for curves of genus $g\ge 2$. 
 
 In what follows we will use the notation of the paper \cite{ap} for simple semihomogeneous bundles. As usual, given a p.p.a.v. $A,\OO_A(\Theta))$ we will denote $\utheta\in NS(A)$ the class of $\OO_A(\Theta)$. Given a rational number $\delta=\frac a b$, with $b>0$, we  denote $E_{A,\delta\utheta}$ a simple semihomogeneous vector bundle on $A$ belonging to Mukai's class $\mathbb S_\delta$, namely such that 
 \[\frac{[\det E]}{\rk E}:=\delta(E)=\delta\utheta\]
  in $NS(A)_\mathbb Q$ (\cite[Theorem 7.11]{semihom}).\footnote{Such simple semihomogeneous vector bundles are not unique,  in the first place because their determinant is not fixed (only its N\'eron-Severi class is). However we will neglect this in the notation, since it is irrelevant in the statement of Theorem \ref{th:variant}} As shown in the paper \cite{ap}, in some sense such bundles can be thought as integral replacements of the $\mathbb Q$-twisted sheaves $\OO_A\langle \lambda L\rangle$. Then we have 
 
 \begin{theorem}\label{th:variant} Let $(A,\OO_A(\Theta))$ be a p.p.a.v., and let  $n\ge 2$.
 The following are equivalent:
\begin{enumerate}
\item $A$ is isomorphic, as polarized variety, to $E\times B$, where $\dim E=1$.
\item For every  $x\in A$ there exists  $\alpha\in \Pic0 A$ such that the graded $\mathcal S_{\OO_A(n\Theta)}$-module 
\begin{equation}\label{eq:module}
\mathcal R_{A,\OO_A(n\Theta)}(\I_x\otimes E_{A,\frac{2n-1} {n-1}\utheta}\otimes P_\alpha)
\end{equation}
 is not generated in degree $0$. 
\item There exists  $x\in A$ and $\alpha\in \Pic0 A$ such that the module \emph{(\ref{eq:module})} is not generated in degree $0$.
\end{enumerate}
\end{theorem}

This is indeed a generalization of Theorem \ref{th:first} because for $n=2$ any simple semihomogeneous vector bundle $E_{A, \frac 3 1\utheta} $ is simply a line bundle algebraically equivalent to $\OO_A(3\Theta)$. To better explain the statement we invoke  the \emph{base point freeness threshold}  (\cite{caucci}), an invariant already mentioned in the introduction and Subsection \ref{sub:p-jets}. We recall that the base point freeness threshold of the polarization $n\utheta$, i.e. $\epsilon_0(n\utheta)$,  is equal to $\frac 1 { n}$ and therefore, concerning the multiplication map surjectivity threshold $\mu_0$ (see Remark \ref{rem:gen}), we have that, thanks to (\ref{eq:thresholds}), 
\[
\mu_0(n\utheta)=\frac1 {n-1}.
\]
 Since a simple semihomogeneous vector bundle $E_{A,\frac 1 {n-1}n\utheta}= E_{A,\frac n {n-1}\utheta}$ behave cohomologically as the $\mathbb Q$-twisted sheaf $\OO_A\langle \frac n {n-1}\utheta\rangle$ (\cite[Proposition 2.1.2]{ap}),  the multiplication maps of global sections 
\begin{equation}\label{eq:ntheta}
H^0(A,\OO_A(n\Theta))\otimes H^0(A,E_{A,\frac n {n-1}\utheta}\otimes P_\alpha )\rightarrow H^0(A, E_{A,\frac n {n-1}\utheta}(n\Theta)\otimes P_\alpha)
\end{equation}
are surjective for general $\alpha\in\Pic0 A$ but not for all. Actually, as we will see below, the dimension of the source of (\ref{eq:ntheta}) is equal to the dimension of the target, namely $n^{2g}$, hence the maps (\ref{eq:ntheta}) are singular for $\alpha$ belonging to a divisor in $\Pic0 A$. In conclusion, one might informally think that, from the present point of view, for $n\ge 3$,  $E_{A,\frac n {n-1}\utheta}$ is to $\OO_A(n\Theta)$ as $\OO_A(2\Theta)$ is to $\OO_A(2\Theta)$.  The simple semihomogeneous bundle $E_{A,\frac{2n-1}{n-1}\utheta}$ is isomorphic to  $E_{A,\frac{n}{n-1}\utheta}(\Theta)$, therefore $E_{A,\frac{2n-1}{n-1}\utheta}$  is to $\OO_A(n\Theta)$ as $\OO_A(3\Theta)$ is to $\OO_A(2\Theta)$. This explains the statement of Theorem \ref{th:variant}.

\begin{proof} \emph{(outline)} The argument follows the lines of the proof of Theorem \ref{th:first}, but there are some complications appearing in the analog of Lemma \ref{lem:generation} (Lemma \ref{lem:generation2} below). In the first place we note that the first reduction holds without change, namely  Theorem \ref{th:variant} is reduced to the following statement, analogous to Theorem \ref{th:reduction}:

\noindent (*)  \emph{If $\dim A\ge 2$ and $\Theta$ is irreducible then the vector bundles $M_{\OO_A(n\Theta)}\otimes
  E_{A, \frac{n+k(n-1)}{n-1}\utheta}\otimes P_\alpha $ are globally generated for all $k\ge 1$ and $\alpha\in\Pic0 A$.  } 
  
  To see the equivalence with  Theorem \ref{th:variant}, we first note that on a decomposable variety, the simple semihomogeneous vector bundle $E_{A,\delta\utheta}$ is the box product of the  bundles $E_{A_i,\delta\utheta_i}$ on the factors $(A_i,\utheta_i)$ (essentially by definition). Therefore the graded module of the statement of Theorem \ref{th:variant} is generated in degree zero if and only if the same is true on each factor. Now (2) holds  on elliptic curves  because for a point $x$ and $p$ on  an elliptic curve $\bar E$ the sheaf $\I_x\otimes E_{\bar E, \frac{2n-1}{n-1}\underline p}=E_{\bar E, \frac{2n-1}{n-1}\underline p}(-x)$ is precisely a simple semihomogeneous vector bundle in the class $\mathbb S_{\frac {2n-1} {n-1}-1}=\mathbb S_{\frac n {n-1}}$ (see (i) below). Therefore the $0$-th degree structure map of the $\mathcal S_{\OO_E(np)}$-module $\mathcal R_{A,\OO_{\bar E}(np)}(\I_x\otimes E_{\bar E,\frac{2n-1} {n-1}\underline p}\otimes P_\alpha)$ is simply a map as (\ref{eq:ntheta}), which is singular for $\alpha$ belonging to a divisor of $\Pic0\bar E$. 
  
  Conversely, the main implication, i.e. (3) $\Rightarrow$ (1), is proved exactly as in Theorem  \ref{th:reduction}. The main facts to keep in mind are: \\
  (i) $E_{A, \frac{n+k(n-1)}{n-1}\utheta}=E_{A,\frac n {n-1}\utheta}(k\Theta)$ (because $\delta(E\otimes F)=\delta(E)+\delta(F)$); \\
  (ii) the vector bundles $E_{A,\lambda\utheta}$ are globally generated as soon as $\lambda>1$ (\cite{p3} after (2-3))\footnote{Here and in what follows beware that in \cite{oprea} and \cite{p3}  is different, namely the (symmetric) simple semihomogeneous vector bundles $E_{A,\frac a b \utheta}$ are denoted $W_{b,a}$. In other words, the role of $a$ and $b$ is exchanged.} ; \\
  (iii)  the vector bundles $M_{n\Theta}\otimes E_{A,\lambda\utheta}$ are IT(0) if and only if $\lambda>\frac{n}{n-1}$. As mentioned above, this is a consequence of the fact that $\mu_0(n\Theta)=\frac 1 {n-1}$. \\
  Plugging this facts in the appropriated places of the argument, the implication (3) $\Rightarrow$ (1) of Theorem \ref{th:variant} follows from (a) above exactly as Theorem \ref{th:first} follows from Theorem \ref{th:reduction}. 
  
  Next, the key point in the proof of (*) is the following generalization of Lemma \ref{lem:generation}

\begin{lemma}\label{lem:generation2} 
 Let $(A,\utheta)$ be an indecomposable p.p.a.v. of dimension $g\ge 2$, and let $n\ge 2$. \begin{enumerate} 
 
 \item\label{it:1} The subvariety $V^0(A,M_{\OO_A(n\Theta)}\otimes E_{A, \frac n{n-1}\utheta})$ (with its reduced structure) is a prime  divisor of $\Pic0 A$ algebraically equivalent to $(n-1)^2n^{2(g-1)}\widehat{\utheta}$.
 
 \item  the vector bundle $M_{\OO_A(n\Theta)}\otimes E_{A, \frac n{n-1}\utheta}$ is \emph{generated by such a prime divisor}, namely the map
\[
 \mathrm{ev}_U:\bigoplus_{\beta\in U}H^0(A, M_{\OO_A(n\Theta)}\otimes E_{A, \frac n{n-1}\utheta}\otimes P_\beta)\otimes P_\beta^{\vee}\rightarrow M_{\OO_A(2\Theta)}\otimes E_{A, \frac n{n-1}\utheta}
\]
 is surjective for all nonempty open subsets $U\subset \Pic0 A$ meeting $V^0(A,M_{\OO_A(n\Theta)}\otimes E_{A, \frac n{n-1}\utheta})$. 
 \end{enumerate}
 \end{lemma} 
 
 The proof that Lemma \ref{lem:generation2} implies (*) is exactly the same of the implication Lemma \ref{lem:generation} $\Rightarrow$ Theorem \ref{th:reduction} (again recalling (i) above) so we do not repeat it here. This concludes the proof of Theorem \ref{th:variant}.
\end{proof}

\begin{proof} (of Lemma \ref{lem:generation2}) Also for this Lemma the line of the argument is the same, but the replacement of the line bundle $\OO_A(2\Theta)$ with the semihomogeneus bundle $E_{A, \frac n{n-1}\utheta}$ is cause of some extra features, reflected by the factor $(n-1)^2$ in item (\ref{it:1}), which is not visible for $n=2$. 

As in Lemma \ref{lem:generation}, applying the Fourier-Mukai transform $\Phi_{\mathcal P}$ to the twisted evaluation map 
\[e:H^0(A,\OO_A(n\Theta))\otimes E_{A,\frac n {n-1}\utheta}\rightarrow E_{A, \frac n{n-1}\utheta}(n\Theta)\]
we get a map of vector bundles
\[
\Phi_{\mathcal P}(e):H^0(A,\OO_A(n\Theta)\otimes \Phi_{\mathcal P}(E_{A,\frac n {n-1}\utheta})\rightarrow \Phi_{\mathcal P}(E_{A, \frac n{n-1}\utheta}(n\Theta))
\]
which fiberwise is the multiplication of global sections
\begin{equation}\label{eq:m-semihom}
H^0(A,\OO_A(n\Theta)\otimes H^0(A,E_{A,\frac n {n-1}\utheta}\otimes P_\alpha)\rightarrow H^0(A,E_{A, \frac n{n-1}\utheta}(n\Theta)\otimes P_\alpha).
\end{equation}
We recall that, for $a,b>0$
\begin{equation}\label{eq:proprieta}
h^0(E_{A,\frac a b\theta})=\chi(E_{A,\frac a b\theta})=a^g, \qquad \mathrm{rk}(E_{A, \frac a b\theta})=b^g, \qquad c_1(E_{A,\frac a b\utheta})=\frac a b\rk(E_{A,\frac a b\utheta})\utheta=ab^{g-1}\utheta
\end{equation}
(\cite[\S 2.1]{oprea} or \cite[\S 1.5]{ap}). Since, by (i) above, $E_{A,\frac n{n-1}\utheta}(n\Theta)\cong E_{A,\frac{n^2}{n-1}\utheta}$, both the source and the target of (\ref{eq:m-semihom}) have dimension $n^{2g}$. Therefore the support of the subvariety $V^0(A, M_{\OO_A(n\Theta)}\otimes E_{A, \frac n{n-1}\utheta})$ coincides with the support of a divisor, namely $D:=\det (\Phi_\cP(e))_0$. Using the third formula in (\ref{eq:proprieta}) and that, for $a,b>0$, 
\[
\Phi_{\cP}(E_{A,\frac a b\utheta})=E_{\widehat A,-\frac b a\hat\utheta}
\]
 (\cite[Proposition 2.2.1]{ap} or \cite[(4-5)]{p3}), after calculations completely similar to those of Lemma \ref{lem:generation} (between (\ref{eq:D})  and (\ref{eq:E})), one gets that
\[
[D]=(n-1)^2n^{2(g-1)}\hat\utheta .
\]
Again following the original argument we consider the pullback of the map $\Phi_{\cP}(e)$ via the isogeny $\varphi_{n\utheta}$, obtaining 
\[
\varphi_{n\utheta}^*\Phi_{\mathcal P}(e):H^0(A,\OO_A(n\Theta)\otimes \varphi_{n\utheta}^*\Phi_{\mathcal P}(E_{A,\frac n {n-1}\utheta})\rightarrow \varphi_{n\utheta}^*\Phi_{\mathcal P}(E_{A, \frac n{n-1}\utheta}(n\Theta)).
\]
 Fiberwise the map $\varphi_{n\utheta}^*\Phi_{\mathcal P}(e)$ is the multiplication map
\begin{equation}\label{eq:pullback}
H^0(A,\OO_A(n\Theta)\otimes H^0(A, t_{(n-1)x}^*E_{A,\frac n {n-1}\utheta})\rightarrow H^0(A,\OO_A(n\Theta)\otimes t_{(n-1)x}^*E_{A,\frac n {n-1}\utheta}).
\end{equation}
This follows from Mukai's theory of semihomogenous vector bundles, for example \cite[Lemma 6.7]{semihom}, asserting that for a semhomogeneous vector bundle $E$ on an abelian variety $A$
\[
t_{ry}^*E\cong E\otimes \varphi_{c_1(E)}(y)
\]
 for all $y\in A$, where $r=\mathrm{rk}\,E$. Since,  for $E=E_{A,\frac n {n-1}\utheta}$, by (\ref{eq:proprieta})
 \[ry=(n-1)(n-1)^{g-1}y, \quad \hbox{and} \quad \varphi_{c_1(E)}(y)=\varphi_{n\utheta}((n-1)y)
 \]  it follows that
for all $x\in A$, 
\[
t^*_{(n-1)x}E_{A,\frac n {n-1}\utheta}\cong E_{A,\frac n {n-1}\utheta}\otimes \varphi_{n\utheta}(x),
\]
 and (\ref{eq:pullback}) follows. 

By the generalization of Kempf's theorem given in \cite[Theorem A]{p3}, the map (\ref{eq:pullback}) is singular if and only if $(n-1)x\in \Theta_a$ for some $a\in A[n]$ ($n$-division points in $A$). It follows that the divisor of zeroes  of the determinant of the map $\varphi_{n\utheta}^*\Phi_\cP(e)$, say $E$, is supported on \ $
\sum_{a\in A[n]}(n-1)_A^*(\Theta_a)=\sum_{a\in A[n]}\bigl((n-1)_A^*(\Theta)\bigr)_{\frac{a}{n-1}}, 
$ 
where $a/(n-1)$ is a (any) element in $(n-1)_{A}^{-1}(a).$ Now, since we have that $[E] = (n-1)^2n^{2g}\utheta = \varphi^*_{n\utheta} [D],$ we conclude that $E = \varphi_{n\utheta}^{\ast}D.$ It follows then that:
\begin{equation}
\label{descomp}
 n^{2g}D = \varphi_{n\utheta\ast}\varphi_{n\utheta}^{\ast} D = \sum_{a\in A[n]}\bigl(\varphi_{n\utheta \ast}(n-1)_{A}^{\ast}\Theta\bigr)_{\varphi_{n\theta}(a/(n-1))}.
\end{equation}
We claim that the divisor $D_{0}:=\varphi_{n\utheta\ast}((n-1)_A^*\Theta)$ is prime (i.e reduced and irreducible). To prove our claim we note that, since $\Theta$ is normal (\cite[Theorem 1]{el}) and ample, the same holds for the \'etale cover $(n-1)_A^*\Theta$. Therefore $(n-1)_A^*\Theta$ is normal and connected, hence a prime divisor. Now, since the restriction of the isogeny $\varphi_{n\utheta}$ to $(n-1)_A^*\Theta$ is birational onto its image (for example, the argument in \cite[p.1591]{p3} applies also to this case), also  $D_{0}$ is a prime divisor, which concludes the proof of the claim. Obviously, the same is true for all the translates appearing in \eqref{descomp} and hence the right hand side of \eqref{descomp} is exactly the decomposition of the divisor $n^{2g}D$ in prime divisors. Since there are exactly $n^{2g}$ summands, it follows that all such summands are equal to $D$ and $D$ is prime, as we wanted to prove. 

\end{proof}

%multiplication map of global sections
%\begin{equation}\label{eq:fract}
%H^0(A,\OO_A(n\Theta))\otimes H^0(A,\I_x\otimes E_{A,\frac{2n-1} {n-1}\utheta}\otimes P_\alpha)\rightarrow H^0(A,\I_x(n\Theta)\otimes E_{A,\frac{2n-1} {n-1}\utheta}\otimes P_\alpha)
%\end{equation}

 %%%%%%%%%%%%%
 \subsection{Variants of Theorem \ref{th:second}}\label{sub:variants2} Passing to the second characterization,  as it happens for hyperelliptic curves, there are many other different second order gaussian maps whose non-surjectivity is equivalent to the non-surjectivity of (\ref{eq:g2(alpha)}). Thus their non-surjectivity characterizes as well decomposable abelian varieties of the form $E\times B$.  Rather than giving a general formula (see however \cite[\S6]{nelson}) we show a couple of examples.

\noindent \emph{ (a) A polarized abelian variety $(A,\OO_A(\Theta))$ is isomorphic to $E\times B$ as in Theorem \ref{th:second} if and only if there exists an $\alpha\in\Pic0 A$ such that the second order gaussian map
  \begin{equation}\label{eq:ex1}
  \mathrm{Rel}^2_A(\OO_A(4\Theta),\OO_A(12\Theta)\otimes P_\alpha)\rightarrow H^0(A, S^2\Omega^1_A\otimes \OO_A(16\Theta)\otimes P_\alpha)
 \end{equation}
is not surjective. }

 This follows from the formula (\ref{eq:thresholds}), coupled with Theorem \ref{th:third}. Indeed, starting with the polarization  $\l=4\utheta$, from Theorem \ref{th:third} we have that if $A$ is not of the form $E\times B$ then $\epsilon_2(4\utheta)={\frac 3 4 } \epsilon_2(3\utheta)<\frac 3 4$. By (\ref{eq:thresholds}) this implies that $\mu_2(4\utheta)< 3$ which means that the gaussian maps (\ref{eq:ex1})
are surjective for all $\alpha\in \Pic0 A$. Conversely, it can be shown that for elliptic curves this does not happen, hence it does not happen also for decomposable varieties $A\cong E\times B$.
 
  One has similar characterizations even when $\mu_2(k\utheta)$ is not integer. These involve simple semihomogeneous vector bundles, as in the previous subsections. An example of this is as follows:

  \emph{ (b) A polarized abelian variety $(A,\OO_A(\Theta))$ is isomorphic to $E\times B$ as in Theorem \ref{th:second} if and only if there exists an $\alpha\in\Pic0 A$ such that the second order gaussian maps
  \begin{equation}\label{eq:ex2}
\mathrm{Rel}^2_A(\OO_A(5\Theta),E_{A,{\frac {15} 2}\utheta}\otimes P_\alpha)\rightarrow H^0(A, S^2\Omega^1_A\otimes \OO_A(5\Theta)\otimes E_{A,{\frac {15} 2}\utheta}\otimes P_\alpha)
 \end{equation}
is not surjective. }  
  
  Indeed let  $\l=5\utheta$. Then, as above, if $A$ is not of the form $E\times B$ then $\epsilon_2(5\utheta)<\frac 3 5$. Thus, by (\ref{eq:thresholds}), $\mu_2(5\theta)<\frac 3 2$, which means that: \emph{(*) $R^2_{A,\OO_A(5\Theta)}\langle{\frac{15} 2}\utheta\rangle$ is an IT(0) $\mathbb Q$-twisted coherent sheaf}. Using the method introduced in \cite{ap} this can be expressed as the surjectivity of certain (second order) gaussian maps by means of semihomogeneous vector bundles.  Indeed, considering a simple semihomogeneous vector bundle $E: =E_{A,{\frac {15} 2}\utheta}$, the  condition \emph{(*)} equivalent to the  condition that the locally free sheaf $R^2_{A,\OO_A(5\Theta)}\otimes E_{A,{\frac {15} 2}\utheta}$ is IT(0) (\cite[Proposition 2.1.2]{ap}). In turns this is equivalent to the surjectivity of the gaussian maps (\ref{eq:ex2}). Conversely, it can be shown that for elliptic curves this does not happen, hence it does not happen also for decomposable varieties $A\cong E\times B$.

 \subsection{Relation with the Seshadri constant}\label{sub:sesh} In fact Theorem \ref{th:third} is an effective version of Nakamaye's theorem characterizing polarized  abelian varieties  of the form $E\times B$ as the only ones achieving the minimal Seshadri constant, namely $1$ (\cite[Theorem 1.1]{nakamaye}). 
To explain this, we first note that, in view of (\ref{eq:normalized}), it is natural to consider the \emph{normalized $p$-jets separation thresholds of a polarization $\utheta$}  as the quantities $\frac{\epsilon_p(\utheta)}{p+1}$.
Now let $\epsilon(A,\utheta)$ denote the Seshadri constant of the polarization $\utheta$ (\cite[p.293]{pag1}). 
  %(a)  (\cite[Theorem 3.7]{nelson}) The $p$-jets separation thresholds satisfy the  subadditivity condition: 
 %$\epsilon_{r+p+1}(\utheta)\le \epsilon_r(\utheta)+\epsilon_{p}(\utheta)$.
  Then \cite[Theorem 4.1]{nelson}) states that
 \begin{equation}\label{eq:alvarado}
 \epsilon(A,\utheta)=\sup_{p\ge 0}\frac{p+1}{\epsilon_p(\utheta)}=\lim_{p\rightarrow \infty}\frac{p+1}{\epsilon_p(\utheta)}.
 \end{equation}
Thanks to (\ref{eq:alvarado}), Nakamaye's theorem characterizes polarized abelian varieties $A=E\times B$ as those such that $\sup_{p\ge 0}\frac{p+1}{\epsilon_p(\utheta)}=1$, but in fact Theorem \ref{th:third} chacterizes them as those such that $\frac{3}{\epsilon_2(\utheta)}=1$. In other words, in this case the asymptotic description (\ref{eq:alvarado}) of the Seshadri constant is equivalent to an explicit finite condition, which is turn translated in an explicit determinantal condition via gaussian maps in Theorem \ref{th:second}.  
 
Thus a natural conjecture, already raised in \cite[Question 7.2]{nelson} as a problem, is that something of this sort always happens. Namely for all rational numbers $\lambda\ge 1$ there should be  a positive integer $p_0$ depending only on $g$, $\lambda$ and the type of a given polarization $\utheta$, such that a polarized abelian variety $(A,\utheta)$ is such that  $\epsilon(A,\utheta)\le \lambda$ if and only if $\frac{p_0+1}{\epsilon_{p_0}(\utheta)}\le \lambda$ i.e.  $\epsilon _{p_0}(\theta)\ge \frac{p_0+1}\lambda$, i.e. 
\[\epsilon_{p_0}((p_0+1)\utheta)\ge \frac 1 \lambda.\]  
In turn, via (\ref{eq:thresholds}), this is equivalent to the condition that  
\begin{equation}\label{eq:mu}
\mu_{p_0}((p_0+1)\utheta)\ge \frac 1{\lambda -1}.
\end{equation}
 Thus the above conjecture would imply that the locus of polarized abelian varieties $(A,\OO_A(\utheta))$ such that  
$\epsilon(A,\utheta)\le \lambda$ would be characterized by the non-surjectivity  of suitable gaussian maps.
Again, the complication that such gaussian maps involve some $\mathbb Q$-twisted sheaf can be handled using semihomogeneous vector bundles as in (b) of the previous subsection 
%In fact the condition (\ref{eq:mu}) is equivalent to the fact that, for all $p\le p_0$,  there is some $\alpha\in Pic0 A$ such that the gaussian map of order $p$
%\[
%\mathrm{Rel}^p_A(\OO_A((p+1)\Theta),E_{A,\frac{(p+1}\lambda}{1-\lambda}\otimes P_\alpha)\rightarrow H^0(A, S^{p}\Omega^1_A\otimes \OO_A((p+1)\Theta)\otimes E_{A,\frac{(p+1}\lambda}{1-\lambda}\otimes P_\alpha))
%\]
%is not surjective. 
(of course, in the case concerning hyperelliptic jacobians  mentioned in the introduction, i.e. for $\lambda=2$, this last complication does not appear).

 %%%%%%%%%%%%%%%%

\vskip2truecm\noindent\textbf{Decalarations. }

\noindent 
\textbf{Data availability. } \ This paper has no associated data.

\noindent\textbf{Conflict of interest. } The authors declare that there is no conflict of interest.

 %%%%%%%%%%%%%%%%

\providecommand{\bysame}{\leavevmode\hbox
to3em{\hrulefill}\thinspace}

\end{document}